\documentclass[10pt]{article}
\usepackage[utf8]{inputenc}
\usepackage[english]{babel}
 
\usepackage[dvipsnames]{xcolor}

\usepackage{mleftright}
\usepackage{tensor}
\usepackage[nottoc]{tocbibind}
\usepackage{amsmath}
\usepackage{cases}
\usepackage{amsfonts}
\usepackage{amssymb}
\usepackage{array}
\usepackage{mathtools}
\usepackage{tikz}
\usepackage{tikz-cd}
\usepackage{ragged2e}
\usepackage{amsthm}
\usepackage{tabularx}
\usepackage{centernot}
\usepackage{faktor}
\usepackage{rotating}
\usepackage{commath}
\usepackage{comment}
\usepackage{mathabx}
\usepackage{fge}
\usepackage[parfill]{parskip}
\usepackage[shortlabels]{enumitem}
\usepackage{xifthen}
\usepackage{etoolbox}
\usepackage{verbatim}
\usepackage{tikzducks}
\usepackage{ytableau}
\usepackage{subfig}
\usepackage{slashed}
\usepackage{xfrac}
\usepackage[title]{appendix}
\usepackage{multirow}
\usepackage{textcomp}
\usepackage[hyphens]{url}

\usepackage[margin=1.15in]{geometry}

\usepackage{pdflscape}
\usepackage{afterpage}
\usepackage{capt-of}

\usepackage{lipsum}
\usepackage{hyperref}

\hypersetup{
    colorlinks=true,
    linkcolor=BrickRed,
    filecolor=Apricot,      
    urlcolor=Rhodamine,
    citecolor=ForestGreen,
    anchorcolor=Fuchsia,
}

\DeclareFontFamily{U}{MnSymbolC}{}
\DeclareSymbolFont{MnSyC}{U}{MnSymbolC}{m}{n}
\DeclareFontShape{U}{MnSymbolC}{m}{n}{
    <-6>  MnSymbolC5
   <6-7>  MnSymbolC6
   <7-8>  MnSymbolC7
   <8-9>  MnSymbolC8
   <9-10> MnSymbolC9
  <10-12> MnSymbolC10
  <12->   MnSymbolC12}{}
\DeclareMathSymbol{\hook}{\mathbin}{MnSyC}{'270}

\newcommand{\eps}{
    \epsilon
}
\newcommand{\rarr}{
    \rightarrow
}

\newcommand{\inner}[2]{
    \langle #1,#2 \rangle
}
\newcommand{\binner}[2]{
    \Big\langle #1,#2 \Big\rangle
}

\renewcommand{\del}[1]{
    \partial #1
}
\newcommand{\delbar}[1]{
    \overline{\partial} #1
}

\newcommand{\tr}[1]{
    \text{\normalfont tr}\,#1
}

\newcommand{\wtilde}[1]{
    \widetilde{#1}
}
\newcommand{\what}[1]{
    \widehat{#1}
}

\newcommand{\Rm}{
    \text{\normalfont{Rm}}
}

\renewcommand{\Re}{
    \text{\normalfont{Re}\,}
}

\newtheorem{innercustomgeneric}{\customgenericname}
\providecommand{\customgenericname}{}
\newcommand{\newcustomtheorem}[2]{%
  \newenvironment{#1}[1]
  {%
   \renewcommand\customgenericname{#2}%
   \renewcommand\theinnercustomgeneric{##1}%
   \innercustomgeneric
  }
  {\endinnercustomgeneric}
}

\makeatletter
\providecommand{\subtitle}[1]{
  \apptocmd{\@title}{\par {\large #1 \par}}{}{}
}
\makeatother

\makeatletter
\def\mathcolor#1#{\@mathcolor{#1}}
\def\@mathcolor#1#2#3{%
  \protect\leavevmode
  \begingroup
    \color#1{#2}#3%
  \endgroup
}
\makeatother

\newcommand\numberthis{\addtocounter{equation}{1}\tag{\theequation}}

\newcustomtheorem{customthm}{Theorem}
\newcustomtheorem{customlem}{Lemma}
\newcustomtheorem{customprop}{Proposition}

\theoremstyle{plain}
\newtheorem{thm}{Theorem}[section]

\newtheorem{propn}[thm]{Proposition}
\newtheorem{cor}[thm]{Corollary}

\theoremstyle{definition}

\newtheorem{rmk}[thm]{Remark}

\theoremstyle{remark}

\newcommand\xqed[1]{%
  \leavevmode\unskip\penalty9999 \hbox{}\nobreak\hfill
  \quad\hbox{#1}}
\newcommand\blktr{\xqed{$\blacktriangle$}}

\title{Anomaly Flow: Shi-Type Estimates and Long-time Existence}
\author{Caleb Suan \\ \texttt{calebkw@math.ubc.ca}}
\date{}

\numberwithin{equation}{section}

\begin{document}

\maketitle

\begin{abstract}
\label{abstract}
    We consider the long-time existence of the anomaly flow on a compact complex $3$-fold with general slope parameter $\alpha'$. In particular, we obtain integral Shi-type estimates for the flow by adapting a integration-by-parts type argument instead of the usual maximum principle techniques. Following this, we prescribe a sufficient smallness condition on $\alpha'$ in order to extend the flow on $[0,\tau)$ to $[0,\tau + \eps)$.
\end{abstract}

\tableofcontents

\section{Introduction}
\label{sect-intro}

The Hull--Strominger system \cite{Hul86,Str86} is a system of partial differential equations generalizing the compactification of the $10$-dimensional heterotic string proposed by Candelas--Horowitz--Strominger--Witten \cite{CHSW85}. On a compact complex $3$-fold $X$ with nowhere-vanishing holomorphic $(3,0)$-form $\Omega$ and a holomorphic vector bundle $E \rarr X$, the equations seek Hermitian metrics $\omega$ on $X$ and $H$ on $E$ satisfying
\begin{equation}
\label{eqn-HS-1}
    F^{2,0} = F^{0,2} = 0, \qquad \omega \wedge F^{1,1} = 0,
\end{equation}
\begin{equation}
\label{eqn-HS-2}
    \sqrt{-1} \del \delbar \omega - \alpha' \Big( \tr (\Rm \wedge \Rm) - \tr (F \wedge F) \Big) = 0,
\end{equation}
\begin{equation}
\label{eqn-HS-3}
    d^\dagger \omega = \sqrt{-1} (\delbar - \del) \log \|\Omega\|_\omega,
\end{equation}
where $\Rm$ and $F$ are respectively the Chern curvatures of $\omega$ and $H$. 

The first condition requires that $H$ is a Hermitian Yang--Mills metric with respect to $\omega$. The second is called the heterotic Bianchi identity which comes from the Green--Schwarz anomaly cancellation in string theory \cite{GS84}. Here $\alpha'$ is a constant called the slope parameter. In \cite{LY05}, Li--Yau show that the third condition is equivalent to the metric $\omega$ being conformally balanced, that is
\begin{equation}
\label{eqn-HS-3'}
    d (\|\Omega\|_\omega \omega^2) = 0.
\end{equation}

\begin{rmk}
\label{rmk-cnxns}
    We note that some authors may use other connections instead of the Chern connection when defining the Hull--Strominger system and also the anomaly flow below. The Hull--Strominger system and anomaly flow may also be generalized to other dimensions \cite{PPZ19a} and to manifolds with other special structures \cite{CGFT22,dlOLS18,dSGFLSE24,FIUV15,II05}. Our focus, however, will be in the complex dimension $3$ setting.
    \blktr
\end{rmk}

In an effort to find solutions to the Hull--Strominger system, Phong--Picard--Zhang \cite{PPZ18c} propose the anomaly flow -- a geometric flow of Hermitian metrics given by
\begin{equation}
\label{eqn-anomaly-flow}
    \del_t (\|\Omega\|_\omega \omega^2) = \sqrt{-1} \del \delbar \omega - \alpha' \Big( \tr (\Rm \wedge \Rm) - \Phi \Big)
\end{equation}
where $\Phi$ is a given closed $(2,2)$-form in $c_2(X)$ that evolves with time. By Chern--Weil theory, the flow preserves the conformally balanced condition.

\begin{rmk}
\label{rmk-other-non-Kahler-flows}
    The conformally balanced condition \eqref{eqn-HS-3'} broadens our scope and allows us to consider non-K\"{a}hler Hermitian metrics. Other flows in non-K\"{a}hler complex geometry have been studied as well, such as the pluriclosed flow \cite{GFGMS24,GFJS23,ST10} and Hermitian curvature flow \cite{ST11,Ust18}.
    \blktr
\end{rmk}

\begin{rmk}
\label{rmk-coupled-flow}
    Ideally one would also have a flow of Hermitian metrics $H$ on $E \rarr X$ to accompany the anomaly flow in order to solve the full Hull--Strominger system. For example, as proposed in \cite{PPZ18c}, we could define a flow
    \begin{equation}
        H^{-1} \del_t H = - \Lambda F
    \end{equation}
    and set $\Phi = (\tr F \wedge F)$. It follows that stationary points of the coupled flow satisfy \eqref{eqn-HS-1} - \eqref{eqn-HS-3} given appropriate initial conditions, as desired.
    \blktr
\end{rmk}

Much work has been done on both the Hull--Strominger system and the anomaly flow in various settings (see, for example \cite{AGF12,AMP24,CPY22,FGV21,FHP21b,FP21,FP23,FTY09,GF20,GFGM24,GFRT17,MPS24,OUV17} as well as other works by these authors). In particular, the long-time existence of the anomaly flow has been studied in certain special cases:
\begin{itemize}
    \item when $\alpha' = 0$ over K\"{a}hler manifolds: \cite{PPZ18b},

    \item on $T^2$-fibrations over $K3$ surfaces: \cite{PPZ18a},

    \item on unimodular Lie groups: \cite{PPZ19b},

    \item on $T^4$-fibrations over Riemann surfaces of genus $g \geq 2$: \cite{FHP21a},

    \item on almost-Abelian Lie groups (using non-Chern connections): \cite{Puj21},

    \item on nilmanifolds (using non-Chern connections): \cite{PU21}.
\end{itemize}

The main result of this paper is a long-time existence result of the anomaly flow.
\begin{thm}
    Suppose that there exist positive constants $B$, $C_0$ such that
    \begin{equation}
        B^{-1} \leq \Big( \frac{1}{2 \|\Omega\|_\omega} \Big) \leq B,
    \end{equation}
    \begin{equation}
        |T|, |\overline{T}|, |\Rm|, |DT|, |D \overline{T}| \leq C_0
    \end{equation}
    along the anomaly flow \eqref{eqn-anomaly-flow} on $t \in [0, \tau)$. If $\alpha'$ is sufficiently small (see Theorem \ref{thm-main} for a more precise statement of this condition), then the flow can be extended to $[0, \tau + \eps)$ for some $\eps > 0$.
    \blktr
\end{thm}

The method we use to show this is by obtaining Shi-type estimates of the evolution of covariant derivatives of curvature and torsion as in \cite{PPZ18b}. The added difficulty when $\alpha' > 0$ is that the additional terms have ``too many" derivatives and are not a Laplacian. For example, we have
\begin{equation}
\label{eqn-del_t-Rm-intro}
    \del_t \Rm = \Big( \frac{1}{2 \|\Omega\|_\omega} \Big) \Bigg[ \frac{1}{2} \Delta_R \Rm + \alpha' \Big( \nabla \overline{\nabla} (\Rm * \Rm) \Big) + \text{lower order terms} \Bigg].
\end{equation}
The $\alpha' \Big( \nabla \overline{\nabla} (\Rm * \Rm) \Big)$ term above and other similar ones for higher order estimates are thus unamenable to a maximum principle argument. Such terms also introduce an additional non-linearity in the evolution of $|D^k \Rm|^2$ and $|D^{k+1} T|^2$ for small $k$.

To get around this, we use an integration-by-parts type argument and work with integral norms instead (see \cite{CP22,HW21} for similar arguments). We then apply Gr\"{o}nwall's inequality in place of using the maximum principle to bound $L^{2p}$-norms of covariant derivatives of curvature and torsion.

By appealing to the Sobolev embedding theorem, we can upgrade the $L^p$-estimates to $L^\infty$-estimates. Finally, we apply another standard argument to extend the flow.

\begin{rmk}
\label{rmk-other-Rm*Rm-flows}
    To the author's knowledge, this result is the first to manage an $\alpha' \Big( \nabla \overline{\nabla} (\Rm * \Rm) \Big)$ in \eqref{eqn-del_t-Rm-intro}. A term of this kind arises when the evolution $\del_t g$ of the metric involves an $\alpha' (\Rm * \Rm)$ term. As such, we speculate that this method can be adapted and applied to obtain long-time existence results for flows with a similar issue, such as those found in \cite{CG20,GGI13,MMS24}.
    \blktr
\end{rmk}

\par {\bf Acknowledgements:} The author thanks his supervisor S\'{e}bastien Picard for suggesting this topic of study and also for many helpful conversations and discussions. 

\section{Evolution Equations from the Anomaly Flow}
\label{sect-anomaly-flow}

In this section, we recall the evolution equations of the metric, curvature, and torsion under the anomaly flow. These were first shown by Phong--Picard--Zhang in \cite{PPZ18b}.

\begin{thm}[\cite{PPZ18b} Theorem 1]
\label{thm-del_t-g}
    Under the anomaly flow \eqref{eqn-anomaly-flow}, the metric $g$ evolves by
    \begin{equation}
        \del_t g_{\overline{p} q} = \Big( \frac{1}{2 \|\Omega\|_\omega} \Big) \Big[ - \wtilde{R}_{\overline{p} q}+ g^{\alpha \overline{\beta}} g^{s \overline{r}} T_{\overline{\beta} s q} \overline{T}_{\alpha \overline{r} \overline{p}} - \alpha' g^{s \overline{r}} (\tensor{R}{_[_{\overline{p}}_s^\alpha_\beta} \tensor{R}{_{\overline{r}}_q_]^\beta_\alpha} - \Phi_{\overline{p} s \overline{r} q}) \Big].
    \end{equation}
    Here $\tensor{R}{_{\overline{p}}_q^r_s} = - \del_{\overline{p}} (g^{r \overline{m}} \del_q g_{\overline{m} s})$ is the Chern curvature tensor, $T_{\overline{p} q k} = \del_q g_{\overline{p} k} - \del_k g_{\overline{p} q}$ is the torsion tensor, and $\wtilde{R}_{\overline{p} q} = \tensor{R}{^k_k_{\overline{p}}_q}$ is (one notion of) the Ricci curvature.
    \blktr
\end{thm}

Schematically, we may write this as
\begin{equation}
\label{eqn-del_t-g0}
\begin{aligned}[b]
    \del_t g = \Big( \frac{1}{2 \|\Omega\|_\omega} \Big) \Big[ \Rm + T * \overline{T} + \alpha' \Big( \Rm * \Rm + \Phi \Big) \Big],
\end{aligned}
\end{equation}
where $*$ denotes a finite linear combination of contractions using the metric $g$.

Expanding out the results of Theorems 4 and 5 in \cite{PPZ18b} respectively give the evolution of the curvature and torsion
\begin{equation}
\label{eqn-del_t-Rm0}
\begin{aligned}[b]
    \del_t \Rm &= \Big( \frac{1}{2 \|\Omega\|_\omega} \Big) \Bigg[ \frac{1}{2} \Delta_R \Rm + \nabla \overline{\nabla} (T * \overline{T}) + \overline{\nabla} (T * \Rm) + \nabla (\overline{T} * \Rm) \\
    &\qquad \qquad \qquad + \Rm * \Rm + \overline{\nabla} (T * T * \overline{T}) + \nabla (\overline{T} * \overline{T} * T) \\
    &\qquad \qquad \qquad + T * \overline{T} * \Rm + T * \overline{T} * T * \overline{T} \\
    &\qquad \qquad \qquad + \alpha' \Big( \nabla \overline{\nabla} (\Rm * \Rm) + \nabla \overline{\nabla} \Phi + \Rm * \Phi \\
    &\qquad \qquad \qquad \qquad \qquad + \overline{\nabla} (T * \Rm * \Rm) + \nabla (\overline{T} * \Rm * \Rm) + \Rm * \Rm * \Rm \\
    &\qquad \qquad \qquad \qquad \qquad + \overline{\nabla} (T * \Phi) + \nabla (\overline{T} * \Phi) + T * \overline{T} * \Rm * \Rm + T * \overline{T} * \Phi \Big) \Bigg],
\end{aligned}
\end{equation}
\begin{equation}
\label{eqn-del_t-T0}
\begin{aligned}[b]
    \del_t T &= \Big( \frac{1}{2 \|\Omega|_\omega} \Big) \Bigg[ \frac{1}{2} \Delta_R T + \nabla (T * \overline{T}) + T * \Rm + T * T * \overline{T} \\
    &\qquad \qquad \qquad + \alpha' \Big( \nabla (\Rm * \Rm) + \nabla \Phi + T * \Rm * \Rm + T * \Phi \Big) \Bigg].
\end{aligned}
\end{equation}

For simplicity, we shall write
\begin{equation}
\label{eqn-del_t-Rm}
\begin{aligned}[b]
    \del_t \Rm = \frac{1}{2 \|\Omega\|_\omega} \Bigg[ \frac{1}{2} \Delta_R \Rm + H_1 + \alpha' \Big( \nabla \overline{\nabla} (\Rm * \Rm) + H_2 \Big) \Bigg],
\end{aligned}
\end{equation}
where 
\begin{equation}
\begin{aligned}[b]
    H_1 &= \nabla \overline{\nabla} (T * \overline{T}) + \overline{\nabla} (T * \Rm) + \nabla (\overline{T} * \Rm) + \Rm * \Rm \\
    &\qquad \overline{\nabla} (T * T * \overline{T}) + \nabla (\overline{T} * \overline{T} * T) + T * \overline{T} * \Rm + T * \overline{T} * T * \overline{T},
\end{aligned}
\end{equation}
\begin{equation}
\begin{aligned}[b]
    H_2 &= \nabla \overline{\nabla} \Phi + \Rm * \Phi + \overline{\nabla} (T * \Rm * \Rm) + \nabla (\overline{T} * \Rm * \Rm) + \Rm * \Rm * \Rm \\
    &\qquad + \overline{\nabla} (T * \Phi) + \nabla (\overline{T} * \Phi) + T * \overline{T} * \Rm * \Rm + T * \overline{T} * \Phi.
\end{aligned}
\end{equation}

Similarly, we write
\begin{equation}
\label{eqn-del_t-T}
\begin{aligned}[b]
    \del_t T &= \frac{1}{2 \|\Omega|_\omega} \Bigg[ \frac{1}{2} \Delta_R T + K_1 + \alpha' \Big( \nabla (\Rm * \Rm) + K_2 \Big) \Bigg],
\end{aligned}
\end{equation}
where
\begin{equation}
\begin{aligned}[b]
    K_1 &= \nabla (T * \overline{T}) + T * \Rm + T * T * \overline{T},
\end{aligned}
\end{equation}
\begin{equation}
\begin{aligned}[b]
    K_2 &= \nabla \Phi + T * \Rm * \Rm + T * \Phi.
\end{aligned}
\end{equation}

By the evolution of the Chern connection, we also have expressions for the evolution of covariant derivatives of both curvature and torsion
\begin{equation}
\label{eqn-del_t-nabla^m-nabla-bar^l-Rm}
\begin{aligned}[b]
    \del_t (\nabla^m \overline{\nabla}^l \Rm) &= \sum_{i+j>0} \sum_{i=0}^m \sum_{j=0}^l (\nabla^{m-i} \overline{\nabla}^{l-j} \Rm) * (\nabla^i \overline{\nabla}^j \del_t g) \\
    &\qquad + \nabla^m \overline{\nabla}^l \Bigg( \frac{1}{2 \|\Omega\|_\omega} \Bigg[ \frac{1}{2} \Delta_R \Rm + H_1 + \alpha' \Big( \nabla \overline{\nabla} (\Rm * \Rm) + H_2 \Big) \Bigg] \Bigg),
\end{aligned}
\end{equation}
\begin{equation}
\label{eqn-del_t-nabla^m-nabla-bar^l-T}
\begin{aligned}[b]
    \del_t (\nabla^m \overline{\nabla}^l T) &= \sum_{i+j>0} \sum_{i=0}^m \sum_{j=0}^l (\nabla^{m-i} \overline{\nabla}^{l-j} T) * (\nabla^i \overline{\nabla}^j \del_t g) \\
    &\qquad + \nabla^m \overline{\nabla}^l \Bigg( \frac{1}{2 \|\Omega\|_\omega} \Bigg[ \frac{1}{2} \Delta_R T + K_1 + \alpha' \Big( \nabla (\Rm * \Rm) + K_2 \Big) \Bigg] \Bigg).
\end{aligned}
\end{equation}

These evolution equations will be heavily used in the sequel where we obtain our Shi-type estimates.

\section{Integral Shi-Type Estimates}
\label{sect-shi-type-ests}

Our objective in this section is to obtain uniform $L^\infty$-estimates on covariant derivatives of the curvature and torsion along the anomaly flow under the assumption that the lower derivatives are also uniformly bounded.

\subsection{Starting Assumptions}
\label{subsect-assumptions}

Suppose for $k \geq 1$ that there exist positive constants $B, C_0, C_1, \ldots, C_{k-1}$ such that
\begin{equation}
\label{eqn-assumptions-1}
    B^{-1} \leq \Big( \frac{1}{2 \|\Omega\|_\omega} \Big) \leq B,
\end{equation}
\begin{equation}
\label{eqn-assumptions-2}
    |D^q \Rm|, |D^{q+1} T|, |D^{q+1} \overline{T}| \leq C_q \text{ for } 1 \leq q \leq k-1,
\end{equation}
\begin{equation}
\label{eqn-assumptions-3}
    |T|, |\overline{T}|, |\Rm|, |DT|, |D \overline{T}| \leq C_0,
\end{equation}
along the anomaly flow on $t \in [0,\tau)$, where
\begin{equation}
    |D^q A|^2 = \sum_{m+l=q} |\nabla^m \overline{\nabla}^l A|^2.
\end{equation}

\begin{rmk}
\label{rmk-dilaton-bounds}
    We note that the first assumption \eqref{eqn-assumptions-1} and the third assumption \eqref{eqn-assumptions-3} in conjunction with \eqref{eqn-del_t-g0} imply that $|\del_t g|$ is uniformly bounded and so the evolving metric $g$ is uniformly bounded above and below by the initial metric. As such, the volume $\int_X 1$ is uniformly bounded along the flow.
    \blktr
\end{rmk}

For simplicity, we shall assume that the evolving form $\Phi$ satisfies any required regularity conditions.

Our goal is to prove that if the assumptions \eqref{eqn-assumptions-1} - \eqref{eqn-assumptions-3} hold for some $k$, then there exists some positive constant $C_k$ such that
\begin{equation}
\label{eqn-C_k}
    |D^k \Rm|, |D^{k+1} T|, |D^{k+1} \overline{T}| \leq C_k.
\end{equation}
As described previously, we do this by first obtaining $L^{2p}$-bounds on $|D^k \Rm|$ and $|D^{k+1} T|$ under appropriate conditions on $\alpha'$. This will be done in \S \ref{subsect-base-ests}. We then take it one step further in \S \ref{subsect-higher-order-ests} and obtain $L^{2p}$-bounds on $|D^{k+1} \Rm|$ and $|D^{k+2} T|$. Finally, we apply the Sobolev embedding theorem and an induction argument to get bounds on the higher order derivatives.

In the sequel, we shall use the convention that $C$ denotes a generic positive constant that may change from line to line but does not depend on the time $t$. The generic constant $C$ {\textbf{MAY}} depend on $\alpha'$, and later $p$ when we start working with $L^{2p}$ norms, and also $\mu$ and $\mu'$ in \S \ref{subsect-k=1}. The constants $B$ and $C_q$ {\textbf{DO NOT}} depend on $\alpha'$, $p$, $\mu$, or $\mu'$. We shall also employ the Einstein summation convention, summing over pairs of matching barred or unbarred upper and lower indices.

By the expression \eqref{eqn-nabla^m-nabla-bar^l-1/2-Omega}, we see that
\begin{equation}
    \Big| D^q \Big( \frac{1}{2 \|\Omega\|_\omega} \Big) \Big| \leq C
\end{equation}
for $0 \leq q \leq k+1$.

We can also check that by \eqref{eqn-nabla-bar^l-nabla^m}
\begin{equation}
    |D^{k+1} \overline{T}| \leq C + |D^{k+1} T|,
\end{equation}
\begin{equation}
\begin{aligned}[b]
    |D^{k+2} \overline{T}| \leq C + C|D^k \Rm| + |D^{k+2} T|,
\end{aligned}
\end{equation}
and so on, hence we only need to worry about $\Rm$ and $T$, and the inequality for $\overline{T}$ will follow in \eqref{eqn-C_k}.

\subsection{Estimates on \texorpdfstring{$|D^k \Rm|$}{|Dk Rm|} and \texorpdfstring{$|D^{k+1} T|$}{|D(k+1) T|}}
\label{subsect-base-ests}

We first obtain pointwise bounds as in \cite{PPZ18b, PPZ18c} before integrating and getting integral bounds. We choose to set $m + l = k \geq 2$ here to ensure that the ensuing terms are no more than quadratic in unknown quantities that we want to bound. The $k = 1$ case will be tackled later in \S \ref{subsect-k=1}.

\subsubsection{Pointwise Estimates}
\label{subsubsect-ptwise-ests}

We compute using \eqref{eqn-del_t-Rm} and \eqref{eqn-del_t-nabla^m-nabla-bar^l-Rm} that
{\allowdisplaybreaks
\begin{align*}
    &\del_t |\nabla^m \overline{\nabla}^l \Rm|^2 \\
    &= C |\del_t g| \cdot |\nabla^m \overline{\nabla}^l \Rm|^2 + 2 \Re \inner{\del_t (\nabla^m \overline{\nabla}^l \Rm)}{\nabla^m \overline{\nabla}^l \Rm} \\
    &\leq C |\del_t g| \cdot |\nabla^m \overline{\nabla}^l \Rm|^2 + \sum_{i+j > 0} \sum_{i=0}^m \sum_{j=0}^l 2 \Re \inner{(\nabla^{m-i} \overline{\nabla}^{l-j} \Rm) * (\nabla^i \overline{\nabla}^j \del_t g)}{\nabla^m \overline{\nabla}^l \Rm} \\*
    &\qquad + 2 \Re \binner{\nabla^m \overline{\nabla}^l \Big( \frac{1}{2 \|\Omega\|_\omega} H_1 \Big)}{\nabla^m \overline{\nabla}^l \Rm} + 2 \alpha' \Re \binner{\nabla^m \overline{\nabla}^l \Big( \frac{1}{2 \|\Omega\|_\omega} H_2 \Big)}{\nabla^m \overline{\nabla}^l \Rm} \\*
    &\qquad + \sum_{i+j < k} \sum_{i=0}^m \sum_{j=0}^l 2 \Re \binner{\Big( \nabla^{m-i} \overline{\nabla}^{l-j} \Big( \frac{1}{2 \|\Omega\|_\omega} \Big) \Big) * \Big( \nabla^i \overline{\nabla}^j \Big( \frac{1}{2} \Delta_R \Rm \Big) \Big)}{\nabla^m \overline{\nabla}^l \Rm} \\*
    &\qquad + \sum_{i+j < k} \sum_{i=0}^m \sum_{j=0}^l 2 \alpha' \Re \binner{\Big( \nabla^{m-i} \overline{\nabla}^{l-j} \Big( \frac{1}{2 \|\Omega\|_\omega} \Big) \Big) * \Big( \nabla^i \overline{\nabla}^j \Big( \nabla \overline{\nabla} (\Rm * \Rm) \Big) \Big)
    }{\nabla^m \overline{\nabla}^l \Rm} \\*
    &\qquad + \Big( \frac{1}{2 \|\Omega\|_\omega} \Big) 2 \Re \binner{\nabla^m \overline{\nabla}^l \Big( \frac{1}{2} \Delta_R \Rm \Big)}{\nabla^m \overline{\nabla}^l \Rm} + \Big( \frac{\alpha'}{2 \|\Omega\|_\omega} \Big) 2 \Re \binner{\nabla^m \overline{\nabla}^l \Big( \nabla \overline{\nabla} (\Rm * \Rm) \Big)}{\nabla^m \overline{\nabla}^l \Rm} \\
    &\leq \underbrace{C |\del_t g| \cdot |\nabla^m \overline{\nabla}^l \Rm|^2}_{\mathbf{(I)}} + \underbrace{C \sum_{i+j > 0} \sum_{i=0}^m \sum_{j=0}^l |\nabla^i \overline{\nabla}^j \del_t g| \cdot |\nabla^m \overline{\nabla}^l \Rm|}_{\mathbf{(II)}} \\*
    &\qquad + \underbrace{C \Big| \nabla^m \overline{\nabla}^l \Big( \frac{1}{2 \|\Omega\|_\omega} H_1 \Big) \Big| \cdot |\nabla^m \overline{\nabla}^l \Rm|}_{\mathbf{(III)}} + \underbrace{C \Big| \nabla^m \overline{\nabla}^l \Big( \frac{1}{2 \|\Omega\|_\omega} H_2 \Big) \Big| \cdot |\nabla^m \overline{\nabla}^l \Rm|}_{\mathbf{(IV)}} \\*
    &\qquad + \underbrace{C \sum_{i+j < k} \sum_{i=0}^m \sum_{j=0}^l |\nabla^i \overline{\nabla}^j (\Delta_R \Rm)| \cdot |\nabla^m \overline{\nabla}^l \Rm|}_{\mathbf{(V)}} + \underbrace{C \sum_{i+j < k} \sum_{i=0}^m \sum_{j=0}^l |\nabla^i \overline{\nabla}^j \nabla \overline{\nabla} (\Rm * \Rm)| \cdot |\nabla^m \overline{\nabla}^l \Rm|}_{\mathbf{(VI)}} \\*
    &\qquad + \underbrace{\Big( \frac{1}{2 \|\Omega\|_\omega} \Big) 2 \Re \binner{\nabla^m \overline{\nabla}^l \Big( \frac{1}{2} \Delta_R \Rm \Big)}{\nabla^m \overline{\nabla}^l \Rm}}_{\mathbf{(VII)}} + \underbrace{\Big( \frac{\alpha'}{2 \|\Omega\|_\omega} \Big) 2 \Re \binner{\nabla^m \overline{\nabla}^l \Big( \nabla \overline{\nabla} (\Rm * \Rm) \Big)}{\nabla^m \overline{\nabla}^l \Rm}}_{\mathbf{(VIII)}}. \numberthis
\end{align*}
}

\paragraph{Terms (I) - (VI)}
\label{para-I-VI-Rm}

The terms here are relatively well-behaved. In short, these all have an appropriate amount of derivatives that can be handled. We recall that
\begin{equation}
\label{eqn-del_t-g-2}
    \del_t g = \frac{1}{2 \|\Omega\|_\omega} \Bigg[ \Rm + T * \overline{T} + \alpha' \Big( \Rm * \Rm + \Phi \Big) \Bigg],
\end{equation}
\begin{equation}
\label{eqn-H_1}
\begin{aligned}[b]
    H_1 &= \nabla \overline{\nabla} (T * \overline{T}) + \overline{\nabla} (T * \Rm) + \nabla (\overline{T} * \Rm) + \Rm * \Rm \\
    &\qquad \overline{\nabla} (T * T * \overline{T}) + \nabla (\overline{T} * \overline{T} * T) + T * \overline{T} * \Rm + T * \overline{T} * T * \overline{T},
\end{aligned}
\end{equation}
\begin{equation}
\label{eqn-H_2}
\begin{aligned}[b]
    H_2 &= \nabla \overline{\nabla} \Phi + \Rm * \Phi + \overline{\nabla} (T * \Rm * \Rm) + \nabla (\overline{T} * \Rm * \Rm) + \Rm * \Rm * \Rm \\
    &\qquad + \overline{\nabla} (T * \Phi) + \nabla (\overline{T} * \Phi) + T * \overline{T} * \Rm * \Rm + T * \overline{T} * \Phi.
\end{aligned}
\end{equation}

Using \eqref{eqn-del_t-g-2} and our assumptions from \S \ref{subsect-assumptions}, we see that
\begin{equation}
\begin{aligned}[b]
    |\nabla^i \overline{\nabla}^j \del_t g| \leq
    \begin{cases}
        C, & i+j \leq k-1 \\
        C + C |D^k \Rm|, & i+j = k.
    \end{cases}
\end{aligned}
\end{equation}

As such, we get that
\begin{equation}
\label{eqn-I-Rm}
\begin{aligned}[b]
    \mathbf{(I)} \leq C |D^k \Rm|^2,
\end{aligned}
\end{equation}
\begin{equation}
\label{eqn-II-Rm}
\begin{aligned}[b]
    \mathbf{(II)} \leq C |D^k \Rm| + C |D^k \Rm|^2.
\end{aligned}
\end{equation}

Equations \eqref{eqn-H_1} and \eqref{eqn-H_2} show that both $H_1$ and $H_2$ have terms with at most $2$ derivatives on $T$ and $\overline{T}$ and at most $1$ on $\Rm$. From this, we see that
\begin{equation}
\begin{aligned}[b]
    \Big| \nabla^m \overline{\nabla}^l \Big( \frac{1}{2 \|\Omega\|_\omega}  H_1 \Big) \Big| \leq C + C |D^k \Rm| + C |D^{k+1} T| + C |D^{k+1} \Rm| + C |D^{k+2} T|,
\end{aligned}
\end{equation}
\begin{equation}
\begin{aligned}[b]
    \Big| \nabla^m \overline{\nabla}^l \Big( \frac{1}{2 \|\Omega\|_\omega}  H_2 \Big) \Big| \leq C + C |D^k \Rm| + C |D^{k+1} T| + C |D^{k+1} \Rm| \\
\end{aligned}
\end{equation}

This gives that
\begin{equation}
\label{eqn-III-Rm}
\begin{aligned}[b]
    \mathbf{(III)} &\leq C |D^k \Rm| + C |D^k \Rm|^2 + C |D^k \Rm| \cdot |D^{k+1} T| \\
    &\qquad + C |D^k \Rm| \cdot |D^{k+1} \Rm| + C |D^k \Rm| \cdot |D^{k+2} T|,
\end{aligned}
\end{equation}
\begin{equation}
\label{eqn-IV-Rm}
\begin{aligned}[b]
    \mathbf{(IV)} &\leq C |D^k \Rm| + C|D^k \Rm|^2 + C |D^k \Rm| \cdot |D^{k+1} T| \\
    &\qquad + C |D^k \Rm| \cdot |D^{k+1} \Rm|.
\end{aligned}
\end{equation}

Similar analysis tells us that
\begin{equation}
\label{eqn-V-Rm}
\begin{aligned}[b]
    \mathbf{(V)} \leq C |D^k \Rm| + C |D^k \Rm|^2 + C |D^k \Rm| \cdot |D^{k+1} \Rm|,
\end{aligned}
\end{equation}
\begin{equation}
\label{eqn-VI-Rm}
\begin{aligned}[b]
    \mathbf{(VI)} \leq C |D^k \Rm| + C |D^k \Rm|^2 + C |D^k \Rm| \cdot |D^{k+1} \Rm|.
\end{aligned}
\end{equation}

\paragraph{Term (VII)}
\label{para-VII-Rm}

The general method for dealing with the term $\mathbf{(VII)}$ is to isolate the highest-order terms (the lower order ones can combine with earlier ones) as well as extract a Laplacian and good negative terms.

The commutator identity \eqref{eqn-nabla^m-nabla-bar^l-Laplace} gives
\begin{equation}
\begin{aligned}[b]
    \nabla^m \overline{\nabla}^l (\Delta_R \Rm) &= \Delta_R (\nabla^m \overline{\nabla}^l \Rm) + \sum_{i=0}^m \sum_{j=0}^l (\nabla^{m-i} \overline{\nabla}^{l-j} \Rm) * (\nabla^i \overline{\nabla}^j \Rm) \\
    &\qquad + \sum_{i=0}^m \sum_{j=0}^l (\nabla^{m-i} \overline{\nabla}^{l+1-j} \Rm) * (\nabla^i \overline{\nabla}^j T) + \sum_{i=0}^m \sum_{j=0}^l (\nabla^{m+1-i} \overline{\nabla}^{l-j} \Rm) * (\nabla^i \overline{\nabla}^j \overline{T}).
\end{aligned}
\end{equation}

This gives
\begin{equation}
\begin{aligned}[b]
    &\Big( \frac{1}{2 \|\Omega\|_\omega} \Big) 2 \Re \binner{\nabla^m \overline{\nabla}^l \Big( \frac{1}{2} \Delta_R \Rm \Big)}{\nabla^m \overline{\nabla}^l \Rm} \\
    &\leq \Big( \frac{1}{2 \|\Omega\|_\omega} \Big) \Re \inner{\Delta_R (\nabla^m \overline{\nabla}^l \Rm)}{\nabla^m \overline{\nabla}^l \Rm} \\
    &\qquad + C \sum_{i=0}^m \sum_{j=0}^l |\nabla^{m-i} \overline{\nabla}^{l-j} \Rm| \cdot |\nabla^i \overline{\nabla}^j \Rm| \cdot |\nabla^m \overline{\nabla}^l \Rm| \\
    &\qquad + C \sum_{i=0}^m \sum_{j=0}^l |\nabla^{m-i} \overline{\nabla}^{l+1-j} \Rm| \cdot |\nabla^i \overline{\nabla}^j T| \cdot |\nabla^m \overline{\nabla}^l \Rm| \\
    &\qquad + C \sum_{i=0}^m \sum_{j=0}^l |\nabla^{m+1-i} \overline{\nabla}^{l-j} \Rm| \cdot |\nabla^i \overline{\nabla}^j \overline{T}| \cdot |\nabla^m \overline{\nabla}^l \Rm| \\
    &\leq \Big( \frac{1}{2 \|\Omega\|_\omega} \Big) \Re \inner{\Delta_R (\nabla^m \overline{\nabla}^l \Rm)}{\nabla^m \overline{\nabla}^l \Rm} \\
    &\qquad + C |D^k \Rm| + C |D^k \Rm|^2 + C |D^k \Rm| \cdot |D^{k+1} \Rm|.
\end{aligned}
\end{equation}

One can then check that (since $m \leq k$)
\begin{equation}
\begin{aligned}[b]
    &2 \Re \inner{\Delta_R (\nabla^m \overline{\nabla}^l \Rm)}{\nabla^m \overline{\nabla}^l \Rm} \\
    &= \inner{\Delta_R (\nabla^m \overline{\nabla}^l \Rm)}{\nabla^m \overline{\nabla}^l \Rm} + \inner{\nabla^m \overline{\nabla}^l \Rm}{\Delta_R (\nabla^m \overline{\nabla}^l \Rm)} \\
    &= \Delta_R |\nabla^m \overline{\nabla}^l \Rm|^2 - 2 |\nabla^{m+1} \overline{\nabla}^l \Rm|^2 - 2 |\nabla^m \overline{\nabla}^{l+1} \Rm|^2 \\
    &\qquad - 2 \Big( |\overline{\nabla} \nabla^m \overline{\nabla}^l \Rm|^2 - |\nabla^m \overline{\nabla}^{l+1} \Rm|^2 \Big) \\
    &\leq \Delta_R |\nabla^m \overline{\nabla}^l \Rm|^2 - 2 |\nabla^{m+1} \overline{\nabla}^l \Rm|^2 - 2 |\nabla^m \overline{\nabla}^{l+1} \Rm|^2 \\
    &\qquad + C \sum_{i=0}^{m-1} |\nabla^m \overline{\nabla}^{l+1} \Rm| \cdot |\nabla^i \Rm| \cdot |\nabla^{m-1-i} \overline{\nabla}^l \Rm| + C \sum_{i=0}^{m-1} |\nabla^i \Rm|^2 \cdot |\nabla^{m-1-i} \overline{\nabla}^l \Rm|^2 \\
    &\leq \Delta_R |\nabla^m \overline{\nabla}^l \Rm|^2 - 2 |\nabla^{m+1} \overline{\nabla}^l \Rm|^2 - 2 |\nabla^m \overline{\nabla}^{l+1} \Rm|^2 + C + C |D^{k+1} \Rm| \\
\end{aligned}
\end{equation}

This means that
\begin{equation}
\label{eqn-VII-Rm}
\begin{aligned}[b]
    \mathbf{(VII)} &\leq \frac{1}{2} \Big( \frac{1}{2 \|\Omega\|_\omega} \Big) \Delta_R \Big( |\nabla^m \overline{\nabla}^l \Rm|^2 \Big) - B^{-1} |\nabla^{m+1} \overline{\nabla}^l \Rm|^2 - B^{-1} |\nabla^m \overline{\nabla}^{l+1} \Rm|^2 \\
    &\qquad + C + C |D^k \Rm| + C |D^{k+1} \Rm| + C |D^k \Rm|^2 + C |D^k \Rm| \cdot |D^{k+1} \Rm|.
\end{aligned}
\end{equation}

\paragraph{Term (VIII)}
\label{para-VIII-Rm}

This final term is one with ``too many" non-Laplacian derivatives. In order to deal with this, we rewrite it in preparation for integration by parts and application of the divergence theorem. We keep track of the constant in front of the terms that are quadratic in the highest order since we will want to cancel them out later.

Using the commutator (Ricci) identity, we have that
\begin{equation}
\begin{aligned}[b]
    \nabla^m \overline{\nabla}^l \Big( \nabla \overline{\nabla} (\Rm * \Rm) \Big) = \overline{\nabla} \nabla^{m+1} \overline{\nabla}^l (\Rm * \Rm) + \sum_{i=0}^m \sum_{j=0}^l \Big( \nabla^{m-i} \overline{\nabla}^{l-j} (\Rm * \Rm) \Big) * (\nabla^i \overline{\nabla}^j \Rm).
\end{aligned}
\end{equation}

From this, we see that
\begin{equation}
\begin{aligned}[b]
    &\Big( \frac{\alpha'}{2 \|\Omega\|_\omega} \Big) 2 \Re \binner{\nabla^m \overline{\nabla}^l \Big( \nabla \overline{\nabla} (\Rm * \Rm) \Big)}{\nabla^m \overline{\nabla}^l \Rm} \\
    &\leq \Big( \frac{\alpha'}{2 \|\Omega\|_\omega} \Big) 2 \Re \inner{\overline{\nabla} \nabla^{m+1} \overline{\nabla}^l (\Rm * \Rm)}{\nabla^m \overline{\nabla}^l \Rm} \\
    &\qquad + C \sum_{i=0}^m \sum_{j=0}^l |\nabla^{m-i} \overline{\nabla}^{l-j} (\Rm * \Rm)| \cdot |\nabla^i \overline{\nabla}^j \Rm| \cdot |\nabla^m \overline{\nabla}^l \Rm| \\
    &\leq \Big( \frac{\alpha'}{2 \|\Omega\|_\omega} \Big) 2 \Re \Big( \overline{\nabla}_{\overline{j}} \inner{\nabla^{m+1} \overline{\nabla}^l (\Rm * \Rm)}{\nabla^m \overline{\nabla}^l \Rm}^{\overline{j}} \Big) \\
    &\qquad - \Big( \frac{\alpha'}{2 \|\Omega\|_\omega} \Big) 2 \Re \inner{\nabla^{m+1} \overline{\nabla}^l (\Rm * \Rm)}{\nabla^{m+1} \overline{\nabla}^l \Rm} \\
    &\qquad + C |D^k \Rm| + C |D^k \Rm|^2 \\
    &\leq \Big( \frac{\alpha'}{2 \|\Omega\|_\omega} \Big) 2 \Re \Big( \overline{\nabla}_{\overline{j}} \inner{\nabla^{m+1} \overline{\nabla}^l (\Rm * \Rm)}{\nabla^m \overline{\nabla}^l \Rm}^{\overline{j}} \Big) \\
    &\qquad + C |D^k \Rm| + C |D^{k+1} \Rm| + C |D^k \Rm|^2 + C |D^k \Rm| \cdot |D^{k+1} \Rm| \\
    &\qquad + 4 a_0 B C_0 \alpha' |\nabla^{m+1} \overline{\nabla}^l \Rm|^2.
\end{aligned}
\end{equation}
where $a_0 > 1$ is a fixed predetermined constant independent of $k$ arising from the linear combinations obscured by the $*$ contraction notation in \eqref{eqn-del_t-Rm} and \eqref{eqn-del_t-T}.

As such, we get
\begin{equation}
\label{eqn-VIII-Rm}
\begin{aligned}[b]
    \mathbf{(VIII)} &\leq 2 \Re \Big( \overline{\nabla}_{\overline{j}} \binner{\Big( \frac{\alpha'}{2 \|\Omega\|_\omega} \Big)  \nabla^{m+1} \overline{\nabla}^l (\Rm * \Rm)}{\nabla^m \overline{\nabla}^l \Rm}^{\overline{j}} \Big) \\
    &\qquad + C |D^k \Rm| + C |D^{k+1} \Rm| + C |D^k \Rm|^2 + C |D^k \Rm| \cdot |D^{k+1} \Rm| \\
    &\qquad + 4 a_0 B C_0 \alpha' |\nabla^{m+1} \overline{\nabla}^l \Rm|^2.
\end{aligned}
\end{equation}

\paragraph{Combining the Terms}
\label{para-combine-Rm}

If we combine what we have from the terms $\mathbf{(I)}$ - $\mathbf{(VIII)}$ (equations \eqref{eqn-I-Rm}, \eqref{eqn-II-Rm}, \eqref{eqn-III-Rm}, \eqref{eqn-IV-Rm}, \eqref{eqn-V-Rm}, \eqref{eqn-VI-Rm}, \eqref{eqn-VII-Rm}, \eqref{eqn-VIII-Rm}), we get the pointwise estimate that
\begin{equation}
\begin{aligned}[b]
    \del_t |\nabla^m \overline{\nabla}^l \Rm|^2 &\leq \frac{1}{2} \Big( \frac{1}{2 \|\Omega\|_\omega} \Big) \Delta_R \Big( |\nabla^m \overline{\nabla}^l \Rm|^2 \Big) - B^{-1} |\nabla^{m+1} \overline{\nabla}^l \Rm|^2 - B^{-1} |\nabla^m \overline{\nabla}^{l+1} \Rm|^2 \\
    &\qquad + 2 \Re \Big( \overline{\nabla}_{\overline{j}} \binner{\Big( \frac{\alpha'}{2 \|\Omega\|_\omega} \Big)  \nabla^{m+1} \overline{\nabla}^l (\Rm * \Rm)}{\nabla^m \overline{\nabla}^l \Rm}^{\overline{j}} \Big) \\
    &\qquad + C + C |D^k \Rm| + C |D^{k+1} \Rm| \\
    &\qquad + C |D^k \Rm|^2 + C |D^k \Rm| \cdot |D^{k+1} T| + C |D^k \Rm| \cdot |D^{k+1} \Rm| + C |D^k \Rm| \cdot |D^{k+2} T| \\
    &\qquad + 4 a_0 B C_0 \alpha' |\nabla^{m+1} \overline{\nabla}^l \Rm|^2.
\end{aligned}
\end{equation}

We have the inequality
\begin{equation}
\begin{aligned}[b]
    |D^{k+1} \Rm|^2 &= \sum_{a+b = k+1} |\nabla^a \overline{\nabla}^b \Rm|^2 \leq \sum_{m+l = k} |\nabla^{m+1} \overline{\nabla}^l \Rm|^2 + |\nabla^m \overline{\nabla}^{l+1} \Rm|^2.
\end{aligned}
\end{equation}

Summing over $m+l = k$ and using the Peter--Paul version of Young's inequality
\begin{equation}
    ab \leq \frac{\eps^{-1}}{2} a^2 + \frac{\eps}{2} b^2 \text{ for } \eps > 0,
\end{equation}
we get for $0 < \eps < 1$ that
{\allowdisplaybreaks
\begin{align*}
\label{eqn-del_t-D^k-Rm}
    \del_t |D^k \Rm|^2 &\leq \frac{1}{2} \Big( \frac{1}{2 \|\Omega\|_\omega} \Big) \Delta_R \Big( |D^k \Rm|^2 \Big)- B^{-1} |D^{k+1} \Rm|^2 \\*
    &\qquad + \sum_{m+l=k} 2 \Re \Big( \overline{\nabla}_{\overline{j}} \binner{\Big( \frac{\alpha'}{2 \|\Omega\|_\omega} \Big) \nabla^{m+1} \overline{\nabla}^l (\Rm * \Rm)}{\nabla^m \overline{\nabla}^l \Rm}^{\overline{j}} \Big) \\*
    &\qquad + C \eps^{-1} \Big( 1 + |D^k \Rm|^2 + |D^{k+1} T|^2 \Big) \\*
    &\qquad + \Big[ C \eps + 4 a_0 B C_0 \alpha' \Big] |D^{k+1} \Rm|^2 + C \eps |D^{k+2} T|^2. \numberthis
\end{align*}
}

A similar treatment for the evolution of covariant derivatives of torsion yields
\begin{equation}
\label{eqn-del_t-D^k+1-T}
\begin{aligned}[b]
    \del_t |D^{k+1} T|^2 &\leq \frac{1}{2} \Big( \frac{1}{2 \|\Omega\|_\omega} \Big) \Delta_R \Big( |D^{k+1} T|^2 \Big) - B^{-1} |D^{k+2} T| \\
    &\qquad + \sum_{m'+l'=k+1} 2 \Re \Big( \nabla_i \binner{\Big( \frac{\alpha'}{2 \|\Omega\|_\omega} \Big) \nabla^{m'} \overline{\nabla}^{l'} (\Rm * \Rm)}{\nabla^{m'} \overline{\nabla}^{l'} T}^i \Big) \\
    &\qquad + C \eps^{-1} \Big( 1 + |D^k \Rm|^2 + |D^{k+1} T|^2 \Big) \\
    &\qquad + \Big[ C \eps + 2 a_0 B C_0 \alpha' \Big] \Big( |D^{k+1} \Rm|^2 + |D^{k+2} T|^2 \Big).
\end{aligned}
\end{equation}

\subsubsection{\texorpdfstring{$L^{2p}$}{L2p}-Estimates}
\label{subsubsect-L^2p-ests-base}

We now define a convenient test function by
\begin{equation}
\label{eqn-G_k}
    G_k = |D^k \Rm|^2 + |D^{k+1} T|^2.
\end{equation}

Summing together \eqref{eqn-del_t-D^k-Rm} and \eqref{eqn-del_t-D^k+1-T}, we get
\begin{equation}
\label{eqn-del_t-G_k}
\begin{aligned}[b]
    \del_t G_k &\leq \frac{1}{2} \Big( \frac{1}{2 \|\Omega\|_\omega} \Big) \Delta_R G_k - B^{-1} G_{k+1} \\
    &\qquad + \sum_{m+l=k} 2 \Re \Big( \overline{\nabla}_{\overline{j}} \binner{\Big( \frac{\alpha'}{2 \|\Omega\|_\omega} \Big) \nabla^{m+1} \overline{\nabla}^l (\Rm * \Rm)}{\nabla^m \overline{\nabla}^l \Rm}^{\overline{j}} \Big) \\
    &\qquad + \sum_{m'+l'=k+1} 2 \Re \Big( \nabla_i \binner{\Big( \frac{\alpha'}{2 \|\Omega\|_\omega} \Big) \nabla^{m'} \overline{\nabla}^{l'} (\Rm * \Rm)}{\nabla^{m'} \overline{\nabla}^{l'} T}^i \Big) \\
    &\qquad + C \eps^{-1} \Big( 1 + G_k \Big) + \Big[ C \eps + 6 a_0 B C_0 \alpha' \Big] G_{k+1}.
\end{aligned}
\end{equation}

We then get that for $p \geq 3$
\begin{equation}
\begin{aligned}[b]
    \del_t \Bigg( \int_X G_k^p \Big) &\leq \int_X \del_t G_k^p + C \int_X G_k^p \\
    &\leq p \int_X G_k^{p-1} \cdot \del_t G_k + C \int_X G_k^p \\
    &\leq \frac{p}{2} \int_X \Big( \frac{1}{2 \|\Omega\|_\omega} \Big) G_k^{p-1} \cdot (\Delta_R G_k) - B^{-1} p \int_X G_k^{p-1} \cdot G_{k+1} \\
    &\qquad + \sum_{m+l=k} 2p\, \Re \Bigg( \int_X G_k^{p-1} \cdot \overline{\nabla}_{\overline{j}} \binner{\Big( \frac{\alpha'}{2 \|\Omega\|_\omega} \Big) \nabla^{m+1} \overline{\nabla}^l (\Rm * \Rm)}{\nabla^m \overline{\nabla}^l \Rm}^{\overline{j}} \Bigg) \\
    &\qquad + \sum_{m'+l'=k+1} 2p\, \Re \Bigg( \int_X G_k^{p-1} \cdot \nabla_i \binner{\Big( \frac{\alpha'}{2 \|\Omega\|_\omega} \Big) \nabla^{m'} \overline{\nabla}^{l'} (\Rm * \Rm)}{\nabla^{m'} \overline{\nabla}^{l'} T}^i \Bigg) \\
    &\qquad + C \eps^{-1} \int_X G_k^{p-1} \cdot \Big( 1 + G_k \Big) + \Big[ C \eps + 6 a_0 B C_0 \alpha' p \Big] \int_X G_k^p \cdot G_{k+1},
\end{aligned}
\end{equation}
where the generic constant $C$ \textbf{MAY} now also depend on $p$. The second term in the first inequality comes from the evolving volume form and that $|\del_t g|$ is bounded along the flow.

\begin{rmk}
\label{rmk-p>=3}
    We impose the condition $p \geq 3$ for now to avoid potentially dividing by $0$ in the future. We cannot just cut out the set $\{G_k = 0\}$ as it would result in a boundary component after integration by parts. The $1 \leq p < 3$ will be addressed later.
    \blktr
\end{rmk}

From
\begin{equation}
    \Delta_R G_k^p = 2 p (p-1) G_k^{p-2} \cdot |\nabla G_k|^2 + p G_k^{p-1} \cdot (\Delta_R G_k),
\end{equation}
we can write the above with an extra negative term.

\begin{equation}
\label{eqn-del_t-int-G_k^p}
\begin{aligned}[b]
    \del_t \Bigg( \int_X G_k^p \Big) &\leq \frac{1}{2} \int_X \Big( \frac{1}{2 \|\Omega\|_\omega} \Big) (\Delta_R G_k^p) - B^{-1} p (p-1) \int_X G_k^{p-2} \cdot |\nabla G_k|^2 - B^{-1} p \int_X G_k^{p-1} \cdot G_{k+1} \\
    &\qquad + \sum_{m+l=k} 2p\, \Re \Bigg( \int_X G_k^{p-1} \cdot \Big( \overline{\nabla}_{\overline{j}} \binner{\Big( \frac{\alpha'}{2 \|\Omega\|_\omega} \Big) \nabla^{m+1} \overline{\nabla}^l (\Rm * \Rm)}{\nabla^m \overline{\nabla}^l \Rm}^{\overline{j}} \Bigg) \\
    &\qquad + \sum_{m'+l'=k+1} 2p\, \Re \Bigg( \int_X G_k^{p-1} \cdot \Big( \nabla_i \binner{\Big( \frac{\alpha'}{2 \|\Omega\|_\omega} \Big) \nabla^{m'} \overline{\nabla}^{l'} (\Rm * \Rm)}{\nabla^{m'} \overline{\nabla}^{l'} T}^i \Bigg) \\
    &\qquad + C \eps^{-1} \int_X G^{p-1} \cdot \Big( 1 + G_k \Big) + \Big[ C \eps + 6 a_0 B C_0 \alpha' p \Big] \int_X G_k^p \cdot G_{k+1}.
\end{aligned}
\end{equation}

Using integration by parts and the divergence theorem \eqref{eqn-div-thm}, we can absorb the Laplacian term and inner product terms into the negative ones.

First, we see that
\begin{equation}
\begin{aligned}[b]
    &\frac{1}{2} \int_X \Big( \frac{1}{2 \|\Omega\|_\omega} \Big) (\Delta_R G_k^p) \\
    &= \frac{1}{2} \int_X \nabla_i \Big[ \Big( \frac{1}{2 \|\Omega\|_\omega} \Big) \cdot g^{i\overline{j}} \cdot \overline{\nabla}_{\overline{j}} G_k^p \Big] - \frac{1}{2} \int_X g^{i\overline{j}} \cdot \nabla_i \Big( \frac{1}{2 \|\Omega\|_\omega} \Big) \cdot \overline{\nabla}_{\overline{j}} G_k^p \\
    &\qquad + \frac{1}{2} \int_X \overline{\nabla}_{\overline{j}} \Big[ \Big( \frac{1}{2 \|\Omega\|_\omega} \Big) \cdot g^{i\overline{j}} \cdot \nabla_i G_k^p \Big] - \frac{1}{2} \int_X g^{i\overline{j}} \cdot \overline{\nabla}_{\overline{j}} \Big( \frac{1}{2 \|\Omega\|_\omega} \Big) \cdot \nabla_i G_k^p \\
    &\leq C \int_X |\nabla G_k^p| \\
    &\leq C \int_X G_k^{p-1} \cdot |\nabla G_k| \\
    &\leq C \eps^{-1} \int_X G_k^{p-2} \cdot G_k^2 + C \eps \int_X G_k^{p-2} \cdot |\nabla G_k|^2 \\
    &= C \eps^{-1} \int_X G_k^p + C \eps \int_X G_k^{p-2} \cdot |\nabla G_k|^2.
\end{aligned}
\end{equation}

Using the divergence theorem \eqref{eqn-div-thm} and that $|D^k \Rm|^2 \leq G_k$, we can write
{\allowdisplaybreaks
\begin{align*}
    &\sum_{m+l=k} 2p\, \Re \Bigg( \int_X G_k^{p-1} \cdot \Big( \overline{\nabla}_{\overline{j}} \binner{\Big( \frac{\alpha'}{2 \|\Omega\|_\omega} \Big) \nabla^{m+1} \overline{\nabla}^l (\Rm * \Rm)}{\nabla^m \overline{\nabla}^l \Rm}^{\overline{j}} \Bigg) \\
    &= \sum_{m+l=k} 2p\, \Re \Bigg( \int_X \overline{\nabla}_{\overline{j}} \binner{G_k^{p-1} \cdot \Big( \frac{\alpha'}{2 \|\Omega\|_\omega} \Big) \nabla^{m+1} \overline{\nabla}^l (\Rm * \Rm)}{\nabla^m \overline{\nabla}^l \Rm}^{\overline{j}} \Bigg) \\*
    &\qquad - \sum_{m+l = k} 2p\, \Re \Bigg( \int_X \overline{\nabla} G_k^{p-1} \cdot \binner{\Big( \frac{\alpha'}{2 \|\Omega\|_\omega} \Big) \nabla^{m+1} \overline{\nabla}^l (\Rm * \Rm)}{\nabla^m \overline{\nabla}^l \Rm} \Bigg) \\
    &\leq C \int_X G_k^{p-1} \cdot |D^k \Rm| \cdot |D^{k+1} (\Rm * \Rm)| \\*
    &\qquad + \sum_{m+l=k} 2 a_0 B \alpha' p (p-1) \int_X G_k^{p-2} \cdot |\overline{\nabla} G_k| \cdot |\nabla^m \overline{\nabla}^l \Rm| \cdot |\nabla^{m+1} \overline{\nabla}^l (\Rm * \Rm)| \\
    &\leq C \int_X G_k^{p-1} \cdot |D^k \Rm| + C \int_X G_k^{p-1} \cdot |D^k \Rm|^2 + C \int_X G_k^{p-1} \cdot |D^k \Rm| \cdot |D^{k+1} \Rm| \\*
    &\qquad + C \int_X G_k^{p-2} \cdot |\overline{\nabla} G_k| \cdot |D^k \Rm| + C \int_X G_k^{p-2} \cdot |\overline{\nabla} G_k| \cdot |D^k \Rm|^2 \\*
    &\qquad + \sum_{m+l=k} 4 a_0 B C_0 \alpha' p (p-1) \int_X G_k^{p-2} \cdot |\overline{\nabla} G_k| \cdot |\nabla^m \overline{\nabla}^l \Rm| \cdot |\nabla^{m+1} \overline{\nabla}^l \Rm| \\
    &\leq C \eps^{-1} \int_X G_k^{p-1} + C \eps^{-1} \int_X G_k^p \\*
    &\qquad + C \eps \int_X G_k^{p-1} \cdot |D^{k+1} \Rm|^2 + C \eps \int_X G_k^{p-2} \cdot |\overline{\nabla} G_k|^2 \\*
    &\qquad + \sum_{m+l=k} 2 a_0 B C_0 \alpha' p (p-1) \int_X G_k^{p-3} \cdot |\overline{\nabla} G_k|^2 \cdot |\nabla^m \overline{\nabla}^l \Rm|^2 \\*
    &\qquad + \sum_{m+l=k} 2 a_0 B C_0 \alpha' p (p-1) \int_X G_k^{p-3} \cdot G_k^2 \cdot |\nabla^{m+1} \overline{\nabla}^l \Rm|^2 \\
    &\leq C \eps^{-1} \int_X G_k^{p-1} \cdot \Big( 1 + G_k \Big) \\*
    &\qquad + \Big[ C \eps + 2 a_0 B C_0 \alpha' p (p-1) \Big] \int_X G_k^{p-1} \cdot |D^{k+1} \Rm|^2 \\*
    &\qquad + \Big[ C \eps + 2 a_0 B C_0 \alpha' p (p-1) \Big] \int_X G_k^{p-2} \cdot |\overline{\nabla} G_k|^2. \numberthis
\end{align*}
}

Similar methods show that the other inner product term has the same upper bound. Substituting these into \eqref{eqn-del_t-int-G_k^p}, and using that
\begin{equation}
    p G_k^{p-1} \leq (p-1) G_k^p + 1,
\end{equation}
we get
\begin{equation}
\begin{aligned}[b]
    \del_t \Bigg( \int_X G_k^p \Bigg) &\leq C \eps^{-1} \int_X \Big( 1 + G_k^p \Big) \\
    &\qquad + \Big[ C \eps + 4 a_0 B C_0 \alpha' p \Big( p+\frac{1}{2} \Big) - B^{-1} p \Big] \int_X G_k^{p-1} \cdot \Big( |D^{k+1} \Rm|^2 + |D^{k+2} T|^2 \Big) \\
    &\qquad + \Big[ C \eps + 4 a_0 B C_0 \alpha' p (p-1) - B^{-1} p (p-1) \Big] \int_X G_k^{p-2} \cdot |\overline{\nabla} G_k|^2.
\end{aligned}
\end{equation}

As such, if
\begin{equation}
    \alpha' < \frac{1}{4 a_0 B^2 C_0 (p+\frac{1}{2})}
\end{equation}
then by choosing $\eps = \eps(k, \alpha', p)$ carefully, we can absorb corresponding terms into the negative ones. This leaves
\begin{equation}
\begin{aligned}[b]
    \del_t \Bigg( \int_X G_k^p \Bigg) \leq C + C \int_X G_k^p,
\end{aligned}
\end{equation}
where we have used that the volume is bounded along the flow. Using Gr\"{o}nwall's inequality, which states that:
\begin{equation}
    u'(t) \leq \beta(t) \cdot u(t) \text{ on } (a,b) \implies u(t) \leq u(a) \exp{ \Bigg( \int_a^t \beta(s) ds \Bigg) } \text{ on } [a,b),
\end{equation}
and setting $u = 1 + \int_X G_k^p$ and $\beta = C$, we have
\begin{thm}
\label{thm-G_k^p}
    Let $k \geq 2$ and $p \geq 3$. Set $G_k = |D^k \Rm|^2 + |D^{k+1} T|^2$. Given the assumptions \eqref{eqn-assumptions-1} - \eqref{eqn-assumptions-3}, if
    \begin{equation}
        \alpha' < \frac{1}{4 a_0 B^2 C_0 (p+\frac{1}{2})}
    \end{equation}
    then there exists some constant $\Lambda_p = \Lambda_p(k, \alpha') > 0$ such that
    \begin{equation}
    \label{eqn-G_k^p}
        \int_X G_k^p(t) \leq \Bigg( 1 + \int_X G_k^p(0)  \Bigg) e^{\Lambda_p t} < \Bigg( 1 + \int_X G_k^p(0) \Bigg) e^{\Lambda_p \tau}.
    \end{equation}
    That is, $\int_X G_k^p(t)$ is uniformly bounded along the flow.
    
    In particular, after taking a $2p$-th root, we get that both
    \begin{equation}
        \Bigg( \int_X |D^k \Rm(t)|^{2p} \Bigg)^{\frac{1}{2p}} \text{ and } \Bigg( \int_X |D^{k+1} T(t)|^{2p} \Bigg)^{\frac{1}{2p}}
    \end{equation}
    are uniformly bounded along the anomaly flow.
    \blktr
\end{thm}

Using uniform boundedness of the volume, we can retrieve the $L^{2p}$-bounds for $1 \leq p < 3$ via H\"{o}lder's inequality.

\begin{cor}
\label{cor-G_k^p}
    Let $k \geq 2$ and $1 \leq p < 3$. Set $G_k = |D^k \Rm|^2 + |D^{k+1} T|^2$. Given the assumptions \eqref{eqn-assumptions-1} - \eqref{eqn-assumptions-3}, if
    \begin{equation}
        \alpha' < \frac{1}{14 a_0 B^2 C_0},
    \end{equation}
    then
    \begin{equation}
        \int_X G_k^p(t) \text{ and both }
        \Bigg( \int_X |D^k \Rm(t)|^{2p} \Bigg)^{\frac{1}{2p}} \text{ and } \Bigg( \int_X |D^{k+1} T(t)|^{2p} \Bigg)^{\frac{1}{2p}}
    \end{equation}
    are uniformly bounded along the anomaly flow.
    \blktr
\end{cor}

\subsection{Estimates on \texorpdfstring{$|D^{k+1} \Rm|$}{|D(k+1) Rm|} and \texorpdfstring{$|D^{k+2} T|$}{|D(k+2) T|}}
\label{subsect-higher-order-ests}

Now suppose that $k \geq 3$. As in \S \ref{subsect-base-ests}, this is done such that terms are at most quadratic in unknowns. The other cases are deferred to \S \ref{subsect-k=2} and \ref{subsect-k=1}.

Recall the test function $G_k$ defined in \eqref{eqn-G_k}. In an analogous manner to \S \ref{subsect-base-ests}, one can check that under the assumptions \eqref{eqn-assumptions-1} - \eqref{eqn-assumptions-3}, we get the pointwise inequality
\begin{equation}
\label{eqn-del_t-G_k+1}
\begin{aligned}[b]
    \del_t G_{k+1} &\leq + \frac{1}{2} \Big( \frac{1}{2 \|\Omega\|_\omega} \Big) (\Delta_R G_{k+1}) - B^{-1} G_{k+2} \\
    &\qquad + \sum_{m+l = k+1} 2 \Re \Big( \overline{\nabla}_{\overline{j}} \binner{\Big( \frac{\alpha'}{2 \|\Omega\|_\omega} \Big) \nabla^{m+1} \overline{\nabla}^l (\Rm * \Rm)}{\nabla^m \overline{\nabla}^l \Rm}^{\overline{j}} \Big) \\
    &\qquad + \sum_{m' + l'= k+2} 2 \Re \Big( \nabla_i \binner{\Big( \frac{\alpha'}{2 \|\Omega\|_\omega} \Big)  \nabla^{m'} \overline{\nabla}^{l'} (\Rm * \Rm)}{\nabla^{m'} \overline{\nabla}^{l'} T}^i \Big) \\
    &\qquad + C (\eps')^{-1} \Big( 1 + G_k + G_{k+1} \Big) + \Big[ C \eps' + 6 a_0 B C_0 \alpha' \Big] G_{k+2},
\end{aligned}
\end{equation}
for $0 < \eps' < 1$.

The above inequality shows that for $p \geq 3$,
\begin{equation}
\label{del_t-int-G_k+1^p}
\begin{aligned}[b]
    \del_t \Bigg( \int_X G_{k+1}^p \Bigg) &\leq + \frac{1}{2} \int_X \Big( \frac{1}{2 \|\Omega\|_\omega} \Big) (\Delta_R G_{k+1}^p) - B^{-1} p \int_X G_{k+1}^{p-1} \cdot G_{k+2} - B^{-1} p (p-1) \int_X G_{k+1}^{p-2} \cdot |\nabla G_{k+1}|^2 \\
    &\qquad + \sum_{m+l = k+1} 2p\, \Re \Bigg( \int_X G_{k+1}^{p-1} \cdot \overline{\nabla}_{\overline{j}} \binner{\Big( \frac{\alpha'}{2 \|\Omega\|_\omega} \Big) \nabla^m \overline{\nabla}^{l+1} (\Rm * \Rm)}{\nabla^m \overline{\nabla}^l \Rm}^{\overline{j}} \Bigg) \\
    &\qquad + \sum_{m'+l'= k+2} 2p\, \Re \Bigg( \int_X G_{k+1}^{p-1} \cdot \nabla_i \binner{\Big( \frac{\alpha'}{2 \|\Omega\|_\omega} \Big) \nabla^{m'} \overline{\nabla}^{l'} (\Rm * \Rm)}{\nabla^{m'} \overline{\nabla}^{l'} T}^i \Bigg) \\
    &\qquad + C (\eps')^{-1} \int_X G_{k+1}^{p-1} \cdot \Big( 1 + G_k + G_{k+1} \Big) + \Big[ C \eps' + 6 a_0 B C_0 \alpha' p \Big] \int_X G_{k+1}^{p-1} \cdot G_{k+2}.
\end{aligned}
\end{equation}

As in the previous section, the Laplacian term can be absorbed by the other terms
\begin{equation}
\begin{aligned}[b]
    \frac{1}{2} \int_X \Big( \frac{1}{2 \|\Omega\|_\omega} \Big) (\Delta_R G_{k+1}^p) &\leq C (\eps')^{-1} \int_X G_{k+1}^p + C \eps' \int_X G_{k+1}^{p-2} \cdot |\nabla G_{k+1}|^2.
\end{aligned}
\end{equation}

Similar calculations show that the inner product terms together can be bound by
\begin{equation}
\begin{aligned}[b]
    &C (\eps')^{-1} \int_X G_{k+1}^{p-1} \cdot (1 + G_k + G_{k+1} \Big) \\
    &\qquad + \Big[ C \eps' + 4 a_0 B C_0 \alpha' p (p-1) \Big] \int_X G_k^{p-1} \cdot G_{k+2} \\
    &\qquad + \Big[ C \eps' + 4 a_0 B C_0 \alpha' p (p-1) \Big] \int_X G_{k+1}^{p-2} \cdot |\overline{\nabla} G_{k+1}|^2.
\end{aligned}
\end{equation}

In tandem with \eqref{del_t-int-G_k+1^p}, we have
\begin{equation}
\begin{aligned}[b]
    \del_t \Bigg( \int_X G_{k+1}^p \Bigg) &\leq C (\eps')^{-1} \int_X \Big( 1 + G_k^p + G_{k+1}^p \Big) \\
    &\qquad + \Big[ C \eps' + 4 a_0 B C_0 \alpha' p \Big(p+\frac{1}{2}\Big) - B^{-1} p \Big] \int_X G_{k+1}^{p-1} \cdot G_{k+2} \\
    &\qquad + \Big[ C \eps' + 4 a_0 B C_0 \alpha' p (p-1) - B^{-1} p (p-1) \Big] \int_X G_{k+1}^{p-2} \cdot |\overline{\nabla} G_{k+1}|^2,
\end{aligned}
\end{equation}
where we have also used that by Young's inequality
\begin{equation}
    p G_{k+1}^{p-1} \cdot G_k \leq (p-1) G_{k+1}^p + G_k^p, \text{ and } 
    p G_{k+1}^{p-1} \leq (p-1) G_{k+1}^p + 1.
\end{equation}

If $\alpha'$ is sufficiently small, then $\eps' = \eps'(k,\alpha',p)$ can be chosen such that
\begin{equation}
    \del_t \Bigg( \int_X G_{k+1}^p \Bigg) \leq C \int_X 1 + C \int_X G_k^p + C \int_X G_{k+1}^p.
\end{equation}

By our assumptions, the volume is bounded along the flow. Further, by Theorem \ref{thm-G_k^p}, under appropriate conditions on $\alpha'$, so is $\int_X G_k^p$. We conclude the following
\begin{cor}
\label{cor-G_k+1^p}
    Let $k \geq 3$ and $G_{k+1} = |D^{k+1} \Rm|^2 + |D^{k+2} T|^2$. Given the assumptions \eqref{eqn-assumptions-1} - \eqref{eqn-assumptions-3}, if:

    \begin{enumerate}[label=(\alph*)]
        \item
        $p \geq 3$ and
        \begin{equation}
            \alpha' < \frac{1}{4 a_0 B^2 C_0 (p+\frac{1}{2})},
        \end{equation}
        then there exists some constant $\Lambda_p' = \Lambda_p'(k, \alpha') > 0$ such that
        \begin{equation}
        \label{eqn-G_k+1^p}
            \int_X G_{k+1}^p(t) \leq \Bigg( 1 + \int_X G_{k+1}^p(0)  \Bigg) e^{\Lambda_p' t} < \Bigg( 1 + \int_X G_{k+1}^p(0) \Bigg) e^{\Lambda_p' \tau}.
        \end{equation}
        That is, $\int_X G_{k+1}^p(t)$ is uniformly bounded along the flow.
    
        In particular, we get that both
        \begin{equation}
            \Bigg( \int_X |D^{k+1} \Rm(t)|^{2p} \Bigg)^{\frac{1}{2p}} \text{ and } \Bigg( \int_X |D^{k+2} T(t)|^{2p} \Bigg)^{\frac{1}{2p}}
        \end{equation}
        are bounded along the anomaly flow;

        \item
        instead $1 \leq p < 3$ and
        \begin{equation}
            \alpha' < \frac{1}{14 a_0 B^2 C_0},
        \end{equation}
        then
        \begin{equation}
            \int_X G_{k+1}^p(t) \text{ and both }
            \Bigg( \int_X |D^{k+1} \Rm(t)|^{2p} \Bigg)^{\frac{1}{2p}} \text{ and } \Bigg( \int_X |D^{k+2} T(t)|^{2p} \Bigg)^{\frac{1}{2p}}
        \end{equation}
        are uniformly bounded along the anomaly flow.
    \end{enumerate}
    \blktr
\end{cor}

\subsection{Sobolev Embedding and Induction}
\label{subsect-Sobolev-induction}

Recall our assumptions from \S \ref{subsect-assumptions}: for some $k \geq 1$ there exist positive constants $B, C_0, C_1 \ldots, C_{k-1}$ such that along the flow on $t \in [0,\tau)$
\begin{equation}
    B^{-1} \leq \Big( \frac{1}{2 \|\Omega\|_\omega} \Big) \leq B,
\end{equation}
\begin{equation}
    |D^q \Rm|, |D^{q+1} T|, |D^{q+1} \overline{T}| \leq C_q \text{ for } 1 \leq q \leq k-1,
\end{equation}
\begin{equation}
    |T|, |\overline{T}|, |\Rm|, |DT|, |D \overline{T}| \leq C_0.
\end{equation}

Recall also that $a_0$ is a predetermined constant from our expressions \eqref{eqn-del_t-Rm} and \eqref{eqn-del_t-T} for $\del_t \Rm$ and $\del_t T$ that is inherent to the anomaly flow.

Theorem \ref{thm-G_k^p} and Corollaries \ref{cor-G_k^p} and \ref{cor-G_k+1^p} have shown that if $k \geq 3$, $p \geq 3$, and if
\begin{equation}
\label{eqn-alpha'}
    \alpha' < \frac{1}{4 a_0 B^2 C_0 (p + \frac{1}{2})},
\end{equation}
then each of
\begin{equation}
    \|D^k \Rm\|_{L^{2p}(M)}, \|D^{k+1} T\|_{L^{2p}(M)}, \|D^{k+1} \Rm\|_{L^{2p}(M)}, \|D^{k+2} T\|_{L^{2p}(M)}
\end{equation}
are bounded uniformly in $t$ along the anomaly flow.

Suppose this holds for some $p$ with $2p > n = 6$. The Sobolev embedding theorem on $D^k \Rm$ and $D^{k+1} T$ then provides $L^\infty$-bounds on $D^k \Rm$ and $D^{k+1} T$ uniform in $t$ along the flow. (See the paragraphs following Lemma 14.3 and also Lemma 14.4 of \cite{Ham82} for more details). That is, there exists some $C_k$ such that
\begin{equation}
    |D^k \Rm|, |D^{k+1} T|, |D^{k+1} \overline{T}| \leq C_k.
\end{equation}

Importantly, the condition \eqref{eqn-alpha'} does not depend on $k$ and so we can induct on $k$ to obtain $L^\infty$-bounds for $D^q \Rm$ and $D^{q+1} T$ for each $q \geq 3$. We have thus shown
\begin{cor}
\label{cor-L^infty}
    Suppose that the assumptions \eqref{eqn-assumptions-1} - \eqref{eqn-assumptions-3} hold for $k = 3$. If
    \begin{equation}
        \alpha' < \frac{1}{14 a_0 B^2 C_0},
    \end{equation}
    then there exist positive constants $C_q$ for $q \geq 3$ such that
    \begin{equation}
        |D^q \Rm|, |D^{q+1} T|, |D^{q+1} \overline{T}| \leq C_q
    \end{equation}
    along the anomaly flow on $t \in [0,\tau)$.
    \blktr
\end{cor}

\begin{rmk}
\label{rmk-consts-dependency}
    We note that the constants $C_q$ from the previous corollary depend on the initial metric $g_0$, the slope parameter $\alpha'$, and the initial bounds $B, C_0, C_1, C_2$.
    \blktr
\end{rmk}

\section{Lowering the Initial Assumptions}
\label{sect-lowering-base-cases}

We now aim to lower the initial value of $k$ in Corollary \ref{cor-L^infty}. To achieve the $k=2$ case, we only need to establish the higher order estimates of \S \ref{subsect-higher-order-ests}. The $k=1$ case requires more work since the base estimates on $|D \Rm|$ and $|D^2 T|$ also need to be reproven.

\subsection{The \texorpdfstring{$k=2$}{k=2} Case}
\label{subsect-k=2}

The procedure in this case is fairly similar to before as we have the base estimates and the result of Theorem \ref{thm-G_k^p} as a starting point.

Analogous calculations to those in \S \ref{subsect-base-ests} show that
\begin{equation}
\label{eqn-del_t-G_3}
\begin{aligned}[b]
    \del_t G_3 &\leq \frac{1}{2} \Big( \frac{1}{2 \|\Omega\|_\omega} \Big) (\Delta_R G_3) - B^{-1} G_4 \\
    &\qquad + \sum_{m+l = 3} 2 \Re \Big( \overline{\nabla}_{\overline{j}} \binner{\Big( \frac{\alpha'}{2 \|\Omega\|_\omega} \Big) \nabla^{m+1} \overline{\nabla}^l (\Rm * \Rm)}{\nabla^m \overline{\nabla}^l \Rm}^{\overline{j}} \Big) \\
    &\qquad + \sum_{m' + l'= 4} 2 \Re \Big( \nabla_i \binner{\Big( \frac{\alpha'}{2 \|\Omega\|_\omega} \Big)  \nabla^{m'} \overline{\nabla}^{l'} (\Rm * \Rm)}{\nabla^{m'} \overline{\nabla}^{l'} T}^i \Big) \\
    &\qquad + C (\eps'^{-1}) \Big( 1 + G_2 + |D^2 \Rm|^4 + G_3 \Big) + \Big[ C \eps' + 6 a_0 B C_0 \alpha' \Big] G_4.
\end{aligned}
\end{equation}

The main difference between this case and the $k \geq 3$ case is that the term
\begin{equation}
    \inner{D^4 (\Rm * \Rm)}{D^4 \Rm} \leq C |D^4 \Rm|^2 + C |D^3 \Rm| \cdot |D^4 \Rm| + C |D^2 \Rm|^2 \cdot |D^4 \Rm|
\end{equation}
is bounded by a term cubic in the unknowns $|D^2 \Rm|$, $|D^3 \Rm|$, and $|D^4 \Rm|$. After applying Young's inequality, the above can be bounded by
\begin{equation}
\begin{aligned}[b]
    \inner{D^4 (\Rm * \Rm)}{D^4 \Rm} \leq C (\eps')^{-1} |D^2 \Rm|^4 + C (\eps')^{-1} |D^3 \Rm|^2 + C \eps' |D^4 \Rm|^2,
\end{aligned}
\end{equation}
which is how the $|D^2 \Rm|^4$ term appears in 
\eqref{eqn-del_t-G_3}.

Using that $|D^2 \Rm|^4 \leq G_2^2$, we see that for $p \geq 3$
\begin{equation}
\begin{aligned}[b]
    \del_t \Bigg( \int_X G_3^p \Bigg) &\leq \frac{1}{2} \int_X \Big( \frac{1}{2 \|\Omega\|_\omega} \Big) (\Delta_R G_3^p) - B^{-1} p \int_X G_3^{p-1} \cdot G_4 - B^{-1} p (p-1) \int_X G_3^{p-2} \cdot |\nabla G_3|^2 \\
    &\qquad + \sum_{m+l = 3} 2p\, \Re \Bigg( \int_X G_3^{p-1} \cdot \overline{\nabla}_{\overline{j}} \binner{\Big( \frac{\alpha'}{2 \|\Omega\|_\omega} \Big) \nabla^{m+1} \overline{\nabla}^l (\Rm * \Rm)}{\nabla^m \overline{\nabla}^l \Rm}^{\overline{j}} \Bigg) \\
    &\qquad + \sum_{m' + l'= 4} 2p\, \Re \Bigg( \int_X G_3^{p-1} \cdot \nabla_i \binner{\Big( \frac{\alpha'}{2 \|\Omega\|_\omega} \Big)  \nabla^{m'} \overline{\nabla}^{l'} (\Rm * \Rm)}{\nabla^{m'} \overline{\nabla}^{l'} T}^i \Bigg) \\
    &\qquad + C (\eps')^{-1} \int_X G_3^{p-1} \cdot \Big( 1 + G_2 + G_2^2 + G_3 \Big) + \Big[ C \eps' + 6 a_0 B C_0 \alpha' p \Big] \int_X G_3^{p-1} \cdot G_4.
\end{aligned}
\end{equation}

The Laplacian term is again well-behaved for our purposes
\begin{equation}
\begin{aligned}[b]
    \frac{1}{2} \int_X \Big( \frac{1}{2 \|\Omega\|_\omega} \Big) (\Delta_R G_3^p) &\leq C (\eps')^{-1}\int_X G_3^p + C \eps' \int_X G_3^{p-2} \cdot |\nabla G_3|^2.
\end{aligned}
\end{equation}

Further, the inner product terms can be bounded by
\begin{equation}
\begin{aligned}[b]
    &C (\eps')^{-1} \int_X G_3^{p-1} \cdot \Big( 1 + G_2^2 + G_3 \Big) \\
    &\qquad + \Big[ C \eps' + 4 a_0 B C_0 \alpha' p (p-1) \Big] \int_X G_3^{p-1} \cdot |D^4 \Rm|^2 \\
    &\qquad + \Big[ C \eps' + 4 a_0 B C_0 \alpha' p (p-1) \Big] \int_X G_3^{p-2} \cdot |\overline{\nabla} G_3|^2.
\end{aligned}
\end{equation}

As before, we can apply Young's inequality to get that
\begin{equation}
\begin{aligned}[b]
    \del_t \Bigg( \int_X G_3^p \Bigg) &\leq C (\eps')^{-1} \int_X \Big( 1 + G_2^p + G_2^{2p} + G_3^p \Big) \\
    &\qquad + \Big[ C \eps' + 4 a_0 B C_0 \alpha' p \Big( p+\frac{1}{2} \Big) - B^{-1} p \Big] \int_X G_3^{p-1} \cdot G_4 \\
    &\qquad + \Big[ C \eps' + 4 a_0 B C_0 \alpha' p (p-1) - B^{-1} p (p-1) \Big] \int_X G_3^{p-2} \cdot |\nabla G_3|^2.
\end{aligned}
\end{equation}

Taking note that we also need $\int_X G_2^{2p}$ to be uniformly bounded, we then conclude from Theorem \ref{thm-G_k^p} and Corollary \ref{cor-G_k^p} that
\begin{cor}
\label{cor-G_3^p}
    Set $k=2$ and $G_3 = |D^3 \Rm|^2 + |D^{3+1} T|^2$. Given the assumptions \eqref{eqn-assumptions-1} - \eqref{eqn-assumptions-3}, if

    \begin{enumerate}[label=(\alph*)]
        \item 
        $p \geq 3$ and
        \begin{equation}
            \alpha' < \frac{1}{4 a_0 B^2 C_0 (2p+\frac{1}{2})},
        \end{equation}
        then there exists some constant $\Lambda_p' = \Lambda_p'(k, \alpha') > 0$ such that
        \begin{equation}
        \label{eqn-G_3^p}
            \int_X G_3(t) \leq \Bigg( 1 + \int_X G_3(0) \Bigg) e^{\Lambda_p' t} < \Bigg( 1 + \int_X G_3(0) \Bigg) e^{\Lambda_p' \tau}.
        \end{equation}
        In particular, we get that both
        \begin{equation}
            \Bigg( \int_X |D^3 \Rm(t)|^{2p} \Bigg)^{\frac{1}{2p}} \text{ and } \Bigg( \int_X |D^4 T(t)|^{2p} \Bigg)^{\frac{1}{2p}}
        \end{equation}
        are bounded along the anomaly flow.

        \item
        instead $1 < p < 3$ and 
        \begin{equation}
            \alpha' < \frac{1}{26 a_0 B^2 C_0},
        \end{equation}
        then
        \begin{equation}
            \int_X G_3^p(t) \text{ and both }
            \Bigg( \int_X |D^{k+1} \Rm(t)|^{2p} \Bigg)^{\frac{1}{2p}} \text{ and } \Bigg( \int_X |D^{k+2} T(t)|^{2p} \Bigg)^{\frac{1}{2p}}
        \end{equation}
        are uniformly bounded along the anomaly flow.
    \end{enumerate}
    \blktr
\end{cor}

\begin{cor}
\label{cor-L^infty-k=2}
    Suppose that the assumptions \eqref{eqn-assumptions-1} - \eqref{eqn-assumptions-3} hold for $k = 2$. If
    \begin{equation}
        \alpha' < \frac{1}{26 a_0 B^2 C_0},
    \end{equation}
    then there exist positive constants $C_q$ for $q \geq 2$ such that
    \begin{equation}
        |D^q \Rm|, |D^{q+1} T|, |D^{q+1} \overline{T}| \leq C_q
    \end{equation}
    along the anomaly flow on $t \in [0,\tau)$.
    \blktr
\end{cor}

As before, the constants $C_q$ will depend on $g_0$, $\alpha'$, $B$, $C_0$, and $C_1$.

\subsection{The \texorpdfstring{$k=1$}{k=1} Case}
\label{subsect-k=1}

We now aim to show the case when $k=1$. As previously mentioned, this case is more complicated as the base estimates first need to be shown.

\subsubsection{Estimates on \texorpdfstring{$|D^k \Rm|$}{|Dk Rm|} and \texorpdfstring{$|D^{k+1} T|$}{|D(k+1) T|}}
\label{subsubsect-base-ests-k=1}

We can check that
\begin{equation}
\begin{aligned}[b]
    \del_t \Big( |D \Rm|^2 + |D^2 T|^2 \Big) &\leq \frac{1}{2} \Big( \frac{1}{2 \|\Omega\|_\omega} \Big) \Delta_R \Big( |D \Rm|^2 + |D^2 T|^2 \Big) - B^{-1} \Big( |D^2 \Rm|^2 + |D^3 T|^2 \Big) \\
    &\qquad + \sum_{m+l = 1} 2 \Re \Big( \overline{\nabla}_{\overline{j}} \binner{\Big( \frac{\alpha'}{2 \|\Omega\|_\omega} \Big) \nabla^{m+1} \overline{\nabla}^l (\Rm * \Rm)}{\nabla^m \overline{\nabla}^l \Rm}^{\overline{j}} \Big) \\
    &\qquad + \sum_{m' + l'= 2} 2 \Re \Big( \nabla_i \binner{\Big( \frac{\alpha'}{2 \|\Omega\|_\omega} \Big)  \nabla^{m'} \overline{\nabla}^{l'} (\Rm * \Rm)}{\nabla^{m'} \overline{\nabla}^{l'} T}^i \Big) \\
    &\qquad + C \eps^{-1} 1 + C \eps^{-1} \Big( |D \Rm|^2 + |D^2 T|^2 \Big) + \Big[ C \eps + 10 a_0 B \alpha' \Big] |D \Rm|^4 \\
    &\qquad + \Big[ C \eps + 6 a_0 B C_0 \alpha' + 4 a_0 B \alpha' \Big] \Big( |D^2 \Rm|^2 + |D^3 T|^2 \Big).
\end{aligned}
\end{equation}

As expected, the quartic $|D \Rm|^4$ term appears, and so we cannot use the same test function as before. To compensate for this term, we use that
\begin{equation}
\begin{aligned}[b]
    \del_t \Big( |\Rm|^2 + |DT|^2 \Big) &\leq \frac{1}{2} \Big( \frac{1}{2 \|\Omega\|_\omega} \Big) \Delta_R \Big( |\Rm|^2 + |DT|^2 \Big) - B^{-1} \Big( |D \Rm|^2 + |D^2 T|^2 \Big) \\
    &\qquad + 2 \Re \Big( \overline{\nabla}_{\overline{j}} \binner{\Big( \frac{\alpha'}{2 \|\Omega\|_\omega} \Big) \nabla (\Rm * \Rm)}{\Rm}^{\overline{j}} \Big) \\
    &\qquad + \sum_{m' + l'= 1} 2 \Re \Big( \nabla_i \binner{\Big( \frac{\alpha'}{2 \|\Omega\|_\omega} \Big)  \nabla^{m'} \overline{\nabla}^{l'} (\Rm * \Rm)}{\nabla^{m'} \overline{\nabla}^{l'} T}^i \Big) \\
    &\qquad + C \eps^{-1} 1 + \Big[ C \eps + 6 a_0 B C_0 \alpha' \Big] \Big( |D \Rm|^2 + |D^2 T|^2 \Big),
\end{aligned}
\end{equation}
and incorporate it into our test function.

Let $\mu > 0$ be a constant to be determined. We consider a test function of the form
\begin{equation}
    G = \Big[ \alpha' \Big( |\Rm|^2 + |DT|^2 \Big) + \mu \Big] \cdot \Big( |D \Rm|^2 + |D^2 T|^2 \Big) = (\alpha' G_0 + \mu) \cdot G_1.
\end{equation}

The two previous calculations tell us that
{\allowdisplaybreaks
\begin{align*}
    \del_t G &= \del_t \Bigg( \Big[ \alpha' \Big( |\Rm|^2 + |DT|^2 \Big) + \mu \Big] \cdot \Big( |D \Rm|^2 + |D^2 T|^2 \Big) \Bigg) \\
    &= \alpha' \del_t \Big( |\Rm|^2 + |DT|^2 \Big) \cdot \Big( |D \Rm|^2 + |D^2 T|^2 \Big) + \Big[ \alpha' \Big( |\Rm|^2 + |DT|^2 \Big) + \mu \Big] \cdot \del_t \Big( |D \Rm|^2 + |D^2 T|^2 \Big) \\
    &\leq \frac{1}{2} \Big( \frac{1}{2 \|\Omega\|_\omega} \Big) \alpha' \Big( |D \Rm|^2 + |D^2 T|^2 \Big) \cdot \Delta_R \Big( |\Rm|^2 + |D T|^2 \Big) \\*
    &\qquad + \frac{1}{2} \Big( \frac{1}{2 \|\Omega\|_\omega} \Big) \Big[ \alpha' \Big( |\Rm|^2 + |DT|^2 \Big) + \mu \Big] \cdot \Delta_R \Big( |D \Rm|^2 + |D^2 T|^2 \Big) \\
    &\qquad - B^{-1} \alpha' \Big( |D \Rm|^2 + |D^2 T|^2 \Big)^2 - B^{-1} \Big[ \alpha' \Big( |\Rm|^2 + |DT|^2 \Big) + \mu \Big] \cdot \Big( |D^2 \Rm|^2 + |D^3 T|^2 \Big) \\*
    &\qquad + 2 \Re \mathbf{(A)} + 2 \Re \mathbf{(B)} \\
    &\qquad + C \eps^{-1} \alpha' \Big( |D \Rm|^2 + |D^2 T|^2 \Big) \\*
    &\qquad + C \eps^{-1} \Big[ \alpha' \Big( |\Rm|^2 + |DT|^2 \Big) + \mu \Big] + C \eps^{-1} \Big[ \alpha' \Big( |\Rm|^2 + |DT|^2 \Big) + \mu \Big] \cdot \Big( |D \Rm|^2 + |D^2 T|^2 \Big) \\
    &\qquad + \Big[ C \eps + 6 a_0 B C_0 \alpha' \Big] \alpha' \Big( |D \Rm|^2 + |D^2 T|^2 \Big)^2 \\*
    &\qquad + \Big[ C \eps + 10 a_0 B \alpha' \Big] \Big[ \alpha' \Big( |\Rm|^2 + |DT|^2 \Big) + \mu \Big] \cdot \Big( |D \Rm|^2 + |D^2 T|^2 \Big)^2 \\*
    &\qquad + \Big[ C \eps + 6 a_0 B C_0 \alpha' + 4 a_0 B \alpha' \Big] \Big[ \alpha' \Big( |\Rm|^2 + |DT|^2 \Big) + \mu \Big] \cdot \Big( |D^2 \Rm|^2 + |D^3 T|^2 \Big) \numberthis
\end{align*}
}where we have used that $|D \Rm|^4 \leq \Big( |D \Rm|^2 + |D^2 T|^2 \Big)^2$ as necessary. In the above, the terms $\mathbf{(A)}$ and $\mathbf{(B)}$ are given by
\begin{equation}
\begin{aligned}[b]
    \mathbf{(A)} &= \alpha' \Big( |D \Rm|^2 + |D^2 T|^2 \Big) \cdot \Big( \overline{\nabla}_{\overline{j}} \binner{\Big( \frac{\alpha'}{2 \|\Omega\|_\omega} \Big) \nabla (\Rm * \Rm)}{\Rm}^{\overline{j}} \Big) \\
    &\qquad + \alpha' \Big( |D \Rm|^2 + |D^2 T|^2 \Big) \cdot \sum_{m' + l'= 1} \Big( \nabla_i \binner{\Big( \frac{\alpha'}{2 \|\Omega\|_\omega} \Big) \nabla^{m'} \overline{\nabla}^{l'} (\Rm * \Rm)}{\nabla^{m'} \overline{\nabla}^{l'} T}^i \Big),
\end{aligned}
\end{equation}
\begin{equation}
\begin{aligned}[b]
    \mathbf{(B)} &= \Big[ \alpha' \Big( |\Rm|^2 + |DT|^2 \Big) + \mu \Big] \cdot \sum_{m+l = 1} \Big( \overline{\nabla}_{\overline{j}} \binner{\Big( \frac{\alpha'}{2 \|\Omega\|_\omega} \Big) \nabla^{m+1} \overline{\nabla}^l (\Rm * \Rm)}{\nabla^m \overline{\nabla}^l \Rm}^{\overline{j}} \Big) \\
    &\qquad + \Big[ \alpha' \Big( |\Rm|^2 + |DT|^2 \Big) + \mu \Big] \cdot \sum_{m' + l'= 2} \Big( \nabla_i \binner{\Big( \frac{\alpha'}{2 \|\Omega\|_\omega} \Big)  \nabla^{m'} \overline{\nabla}^{l'} (\Rm * \Rm)}{\nabla^{m'} \overline{\nabla}^{l'} T}^i \Big).
\end{aligned}
\end{equation}

By further grouping terms and inflating constants, noting that $\alpha' \Big( |\Rm|^2 + |DT|^2 \Big) + \mu \leq 2C_0^2 + \mu$, we get
\begin{equation}
\label{eqn-aux}
\begin{aligned}[b]
    \del_t G &\leq \frac{1}{2} \Big( \frac{1}{2 \|\Omega\|_\omega} \Big) (\Delta_R G) - B^{-1} \alpha' \Big( |D \Rm|^2 + |D^2 T|^2 \Big)^2 - B^{-1} \mu \Big( |D^2 \Rm|^2 + |D^3 T|^2 \Big) \\
    &\qquad + \Big( \frac{\alpha'}{2 \|\Omega\|_\omega} \Big) 2 \Re \binner{\nabla \Big( |\Rm|^2 + |DT|^2 \Big)}{\nabla \Big( |D \Rm|^2 + |D^2 T|^2 \Big)} \\
    &\qquad + 2 \Re \mathbf{(A)} + 2 \Re \mathbf{(B)} \\
    &\qquad + C \eps^{-1} \Big( 1 + G \Big) \\
    &\qquad + \Big[ C \eps + 20 a_0 B C_0^2 (\alpha')^2 + 6 a_0 B C_0 (\alpha')^2 + 10 a_0 B \alpha' \mu \Big] \Big( |D \Rm|^2 + |D^2 T|^2 \Big)^2 \\
    &\qquad + \Big[ C \eps + 12 a_0 B C_0^3 (\alpha')^2 + 8 a_0 B C_0^2 (\alpha')^2 + 6 a_0 B C_0 \alpha' \mu + 4 a_0 B \alpha' \mu \Big] \Big( |D^2 \Rm|^2 + |D^3 T|^2 \Big).
\end{aligned}
\end{equation}
Note that from now on, the generic constant $C$ \textbf{MAY} depend on $\mu$ as well.

From this, it follows that for $p \geq 3$, we have
\begin{equation}
\begin{aligned}[b]
    \del_t \Bigg( \int_X G^p \Bigg) &\leq \frac{1}{2} \int_X \Big( \frac{1}{2 \|\Omega\|_\omega} \Big) (\Delta_R G^p) - B^{-1} p (p-1) \int_X G^{p-2} \cdot |\nabla G|^2 \\
    &\qquad - B^{-1} \alpha' p \int_X G^{p-1} \cdot \Big( |D \Rm|^2 + |D^2 T|^2 \Big)^2 - B^{-1} \mu p \int_X G^{p-1} \cdot \Big( |D^2 \Rm|^2 + |D^3 T|^2 \Big) \\
    &\qquad + 2p\, \Re \Bigg( \int_X G^{p-1} \cdot \Big( \frac{\alpha'}{2 \|\Omega\|_\omega} \Big) \binner{\nabla \Big( |\Rm|^2 + |DT|^2 \Big)}{\nabla \Big( |D \Rm|^2 + |D^2 T|^2 \Big)} \Bigg) \\
    &\qquad + 2p\, \Re \Bigg( \int_X G^{p-1} \cdot \mathbf{(A)} \Bigg) + 2p\, \Re \Bigg( \int_X G^{p-1} \cdot \mathbf{(B)} \Bigg) \\
    &\qquad + C \eps^{-1} \int_X G^{p-1} \cdot \Big( 1 + G \Big) \\
    &\qquad + \Big[ C \eps + 20 a_0 B C_0^2 (\alpha')^2 p \\
    &\qquad \qquad \qquad + 6 a_0 B C_0 (\alpha')^2 p + 10 a_0 B \alpha' \mu p \Big] \int_X G^{p-1} \cdot \Big( |D \Rm|^2 + |D^2 T|^2 \Big)^2 \\
    &\qquad + \Big[ C \eps + 12 a_0 B C_0^3 (\alpha')^2 p + 8 a_0 B C_0^2 (\alpha')^2 p \\
    &\qquad \qquad + 6 a_0 B C_0 \alpha' \mu p + 4 a_0 B \alpha' \mu p \Big] \int_X G^{p-1} \cdot \Big( |D^2 \Rm|^2 + |D^3 T|^2 \Big).
\end{aligned}
\end{equation}

As seen in the previous sections, the Laplacian can be absorbed into the other terms
\begin{equation}
\begin{aligned}[b]
    \frac{1}{2} \int_X \Big( \frac{1}{2 \|\Omega\|_\omega} \Big) (\Delta_R G^p) \leq C \eps^{-1} \int_X G^p + C \eps \int_X G^{p-2} \cdot |\nabla G|^2.
\end{aligned}
\end{equation}

In preparation for dealing with the upcoming terms, we note that by the commutator (Ricci) identity
\begin{equation}
\begin{aligned}[b]
    \Big| \nabla |D \Rm|^2 \Big| &\leq \Big| \nabla \inner{\nabla \Rm}{\nabla \Rm} \Big| + \Big| \nabla \inner{\overline{\nabla} \Rm}{\overline{\nabla} \Rm} \Big| \\
    &\leq |\nabla \Rm| \cdot |\nabla^2 \Rm| + |\nabla \Rm| \cdot |\overline{\nabla} \nabla \Rm| + |\overline{\nabla} \Rm| \cdot |\nabla \overline{\nabla} \Rm| + |\overline{\nabla} \Rm| \cdot |\overline{\nabla}^2 \Rm| \\
    &\leq C |\Rm| \cdot |D \Rm| + 4 |D \Rm| \cdot |D^2 \Rm|.
\end{aligned}
\end{equation}

Likewise
\begin{equation}
\begin{aligned}[b]
    \Big| \nabla |D^2 T|^2 \Big| &\leq \Big| \nabla \inner{\nabla^2 T}{\nabla^2 T} \Big| + \Big| \nabla \inner{\nabla \overline{\nabla} T}{\nabla \overline{\nabla} T} \Big| + \Big| \nabla \inner{\overline{\nabla}^2 T}{\overline{\nabla}^2 T} \Big| \\
    &\leq C |DT| \cdot |D^2 T| + C |D \Rm| \cdot |D^2 T| + 6 |D^2 T| \cdot |D^3 T|.
\end{aligned}
\end{equation}

We also have
\begin{equation}
    \Big| \nabla |\Rm|^2 \Big| \leq 2 |\Rm| \cdot |D \Rm| \text{ and } \Big| \nabla |DT|^2 \Big| \leq C |\Rm| \cdot |DT| + 4 |DT| \cdot |D^2 T|.
\end{equation}

As before, we keep track of the coefficients of the highest order terms.

The inner product term following the negative terms in \eqref{eqn-aux} is bounded by
{\allowdisplaybreaks
\begin{align*}
    &2p\, \Re \Bigg( \int_X G^{p-1} \cdot \Big( \frac{\alpha'}{2 \|\Omega\|_\omega} \Big) \binner{\nabla \Big( |\Rm|^2 + |DT|^2 \Big)}{\nabla \Big( |D \Rm|^2 + |D^2 T|^2 \Big)} \Bigg) \\
    &\leq 2 B \alpha' p \int_X G^{p-1} \cdot \Big( C |\Rm| \cdot |DT| + 2 |\Rm| \cdot |D \Rm| + 4 |DT| \cdot |D^2 T| \Big) \\*
    &\qquad \qquad \qquad \cdot \Big( C |\Rm| \cdot |D \Rm| + C |DT| \cdot |D^2 T| + C |D \Rm| \cdot |D^2 T| \\*
    &\qquad \qquad \qquad \qquad + 4 |D \Rm| \cdot |D^2 \Rm| + 6 |D^2 T| \cdot |D^3 T| \Big). \numberthis
\end{align*}
}In the above, those terms not involving the generic constant $C$ contribute to the highest order terms. We apply Young's inequality to the other terms to get
{\allowdisplaybreaks
\begin{align*}
    &2p\, \Re \Bigg( \int_X G^{p-1} \cdot \Big( \frac{\alpha'}{2 \|\Omega\|_\omega} \Big) \binner{\nabla \Big( |\Rm|^2 + |DT|^2 \Big)}{\nabla \Big( |D \Rm|^2 + |D^2 T|^2 \Big)} \Bigg) \\
    &\leq C \eps^{-1} \int_X G^{p-1} + C \eps^{-1} \int_X G^p \\*
    &\qquad + C \eps \int_X G^{p-1} \cdot \Big( |D \Rm|^2 + |D^2 T|^2 \Big)^2 + C \eps \int_X G^{p-1} \cdot \Big( |D^2 \Rm|^2 + |D^3 T|^2 \Big) \\*
    &\qquad + 16 B \alpha' p \int_X G^{p-1} \cdot |\Rm| \cdot |D \Rm|^2 \cdot |D^2 \Rm| + 24 B \alpha' p \int_X G^{p-1} \cdot |\Rm| \cdot |D \Rm| \cdot |D^2 T| \cdot |D^3 T| \\*
    &\qquad + 32 B \alpha' p \int_X G^{p-1} \cdot |DT| \cdot |D \Rm| \cdot |D^2 T| \cdot |D^2 \Rm| + 48 B \alpha' p \int_X G^{p-1} \cdot |DT| \cdot |D^2 T|^2 \cdot |D^3 T| \\
    &\leq C \eps^{-1} \int_X G^{p-1} + C \eps^{-1} \int_X G^p \\*
    &\qquad + C \eps \int_X G^{p-1} \cdot \Big( |D \Rm|^2 + |D^2 T|^2 \Big)^2 + C \eps \int_X G^{p-1} \cdot \Big( |D^2 \Rm|^2 + |D^3 T|^2 \Big) \\*
    &\qquad + \frac{1}{8} B^{-1} \alpha' p \int_X G^{p-1} \cdot |D \Rm|^4 + 512 B^3 \alpha' p \int_M G^{p-1} \cdot |\Rm|^2 \cdot |D^2 \Rm|^2 \\*
    &\qquad + \frac{1}{8} B^{-1} \alpha' p \int_X G^{p-1} \cdot |D \Rm|^2 \cdot |D^2 T|^2 + 1152 B^3 \alpha' p \int_X G^{p-1} \cdot |\Rm|^2 \cdot |D^3 T|^2 \\*
    &\qquad + \frac{1}{8} B^{-1} \alpha' p \int_X G^{p-1} \cdot |D \Rm|^2 \cdot |D^2 T|^2 + 2048 B^3 \alpha' p \int_X G^{p-1} \cdot |DT|^2 \cdot |D^2 \Rm|^2 \\*
    &\qquad + \frac{1}{8} B^{-1} \alpha' p \int_X G^{p-1} \cdot |D^2 T|^4 + 4608 B^3 \alpha' p \int_X G^{p-1} \cdot |DT|^2 \cdot |D^3 T|^2 \\
    &\leq C \eps^{-1} \int_X G^{p-1} + C \eps^{-1} \int_X G^p \\*
    &\qquad + \Big[ C \eps + \frac{1}{2} B^{-1} \alpha' p \Big] \int_X G^{p-1} \cdot \Big( |D \Rm|^2 + |D^2 T|^2 \Big)^2 \\*
    &\qquad + \Big[ C \eps + 8320 B^3 C_0^2 \alpha' p \Big] \int_X G^{p-1} \cdot \Big( |D^2 \Rm|^2 + |D^3 T|^2 \Big), \numberthis
\end{align*}
}

It remains to deal with the terms involving $\mathbf{(A)}$ and $\mathbf{(B)}$. For part of the $\mathbf{(A)}$ term, we have
{\allowdisplaybreaks
\begin{align*}
    &2p\, \Re \Bigg( \int_X G^{p-1} \cdot \overline{\nabla}_{\overline{j}} \binner{\Big( \frac{\alpha'}{2 \|\Omega\|_\omega} \Big) \nabla (\Rm * \Rm)}{\Rm}^{\overline{j}} \cdot \alpha' \Big( |D \Rm|^2 + |D^2 T|^2 \Big) \Bigg) \\
    &= 2 (\alpha')^2 p, \Re \Bigg( \int_X \overline{\nabla}_{\overline{j}} \binner{G^{p-1} \cdot \Big( \frac{1}{2 \|\Omega\|_\omega} \Big) \cdot \Big( |D \Rm|^2 + |D^2 T|^2 \Big) \cdot \nabla (\Rm * \Rm)}{\Rm}^{\overline{j}} \Bigg) \\*
    &\qquad - 2 (\alpha')^2 p\, \Re \Bigg( \int_X \overline{\nabla} G^{p-1} \cdot \binner{\Big( \frac{1}{2 \|\Omega\|_\omega} \Big) \cdot \Big( |D \Rm|^2 + |D^2 T|^2 \Big) \cdot \nabla (\Rm * \Rm)}{\Rm} \Bigg) \\*
    &\qquad - 2 (\alpha')^2 p\, \Re \Bigg( \int_X G^{p-1} \cdot \binner{\Big( \frac{1}{2 \|\Omega\|_\omega} \Big) \cdot \nabla (\Rm * \Rm)}{\Rm} \cdot \overline{\nabla} \Big( |D \Rm|^2 + |D^2 T|^2 \Big) \Bigg) \\
    &\leq C \int_X G^{p-1} \cdot \Big( |D \Rm|^2 + |D^2 T|^2 \Big) \cdot |D \Rm| \\*
    &\qquad + 4 a_0 B C_0^2 (\alpha')^2 p (p-1) \int_X G^{p-2} \cdot |\overline{\nabla} G| \cdot \Big( |D \Rm|^2 + |D^2 T|^2 \Big) \cdot |D \Rm| \\*
    &\qquad + 4 a_0 B C_0^2 (\alpha')^2 p \int_X G^{p-1} \cdot |D \Rm| \cdot \Big( C |\Rm| \cdot |D \Rm| + C |DT| \cdot |D^2 T| + C |D \Rm| \cdot |D^2 T| \\*
    &\qquad \qquad \qquad \qquad \qquad \qquad \qquad \qquad \qquad + 4 |D \Rm| \cdot |D^2 \Rm| + 6 |D^2 T| \cdot |D^3 T| \Big). \numberthis
\end{align*}
}We apply Young's inequality to the terms in the middle integral and keep track of the coefficients on the highest order contributions of the third integral, while absorbing the rest into good terms with the generic constant $C$.
{\allowdisplaybreaks
\begin{align*}
    &2p\, \Re \Bigg( \int_X G^{p-1} \cdot \overline{\nabla}_{\overline{j}} \binner{\Big( \frac{\alpha'}{2 \|\Omega\|_\omega} \Big) \nabla (\Rm * \Rm)}{\Rm}^{\overline{j}} \cdot \alpha' \Big( |D \Rm|^2 + |D^2 T|^2 \Big) \Bigg) \\
    &\leq C \eps^{-1} \int_X G^p + C \eps \int_X G^{p-1} \cdot \Big( |D \Rm|^2 + |D^2 T|^2 \Big)^2 \\*
    &\qquad + 2 a_0 B C_0^2 (\alpha')^2 p (p-1) \int_X G^{p-3} \cdot \Big( |D \Rm|^2 + |D^2 T|^2 \Big) \cdot |\overline{\nabla} G|^2 \\*
    &\qquad + 2 a_0 B C_0^2 (\alpha')^2 p (p-1) \int_X G^{p-3} \cdot G^2 \cdot |D \Rm|^2 \cdot \Big( |D \Rm|^2 + |D^2 T|^2 \Big) \\*
    &\qquad + 20 a_0 B C_0^2 (\alpha')^2 p \int_X G^{p-1} \cdot \Big( |D \Rm|^2 + |D^2 T|^2 \Big)^2 + 20 a_0 B C_0^2 (\alpha')^2 p \int_X G^{p-1} \cdot \Big( |D^2 \Rm|^2 + |D^3 T|^2 \Big) \\
    &\leq C \eps^{-1} \int_X G^p + C \eps \int_X G^{p-1} \cdot \Big( |D \Rm|^2 + |D^2 T|^2 \Big)^2 \\*
    &\qquad + 2 a_0 B C_0^2 (\alpha')^2 \mu^{-1} p (p-1) \int_X G^{p-2} \cdot |\overline{\nabla} G|^2 + 2 a_0 B C_0^2 (\alpha')^2 p (p-1) \int_X G^{p-1} \cdot \Big( |D \Rm|^2 + |D^2 T|^2 \Big)^2 \\*
    &\qquad + 20 a_0 B C_0^2 (\alpha')^2 p \int_X G^{p-1} \cdot \Big( |D \Rm|^2 + |D^2 T|^2 \Big)^2 + 20 a_0 B C_0^2 (\alpha')^2 p \int_X G^{p-1} \cdot \Big( |D^2 \Rm|^2 + |D^3 T|^2 \Big) \\
    &\leq C \eps^{-1} \int_X G^p + 2 a_0 B C_0^2 (\alpha')^2 \mu^{-1} p (p-1) \int_X G^{p-2} \cdot |\overline{\nabla} G|^2 \\*
    &\qquad + \Big[ C \eps + 2 a_0 B C_0^2 (\alpha')^2 p (p-1) + 20 a_0 B C_0^2 (\alpha')^2 p \Big] \int_X G^{p-1} \cdot \Big( |D \Rm|^2 + |D^2 T|^2 \Big)^2 \\*
    &\qquad + 20 a_0 B C_0^2 (\alpha')^2 p \int_X G^{p-1} \cdot \Big( |D^2 \Rm|^2 + |D^3 T|^2 \Big). \numberthis
\end{align*}
}

The other parts of the $\mathbf{(A)}$ term have similar bounds and so
\begin{equation}
\begin{aligned}[b]
    &2p\, \Re \Bigg( \int_X G^{p-1} \cdot \mathbf{(A)} \Bigg) \\
    &\leq C \eps^{-1} \int_X G^p + 6 B C_0^2 L (\alpha')^2 \mu^{-1} p (p-1) \int_X G^{p-2} \cdot |\overline{\nabla} G|^2 \\
    &\qquad + \Big[ C \eps + 6 B C_0^2 L (\alpha')^2 p (p-1) + 60 B C_0^2 L (\alpha')^2 p \Big] \int_X G^{p-1} \cdot \Big( |D \Rm|^2 + |D^2 T|^2 \Big)^2 \\
    &\qquad + 60 B C_0^2 L (\alpha')^2 p \int_X G^{p-1} \cdot \Big( |D^2 \Rm|^2 + |D^3 T|^2 \Big).
\end{aligned}
\end{equation}

For the $\mathbf{(B)}$ term, we check that
{\allowdisplaybreaks
\begin{align*}
    &2p\, \Re \Bigg( \int_X G^{p-1} \cdot \overline{\nabla}_{\overline{j}} \binner{\Big( \frac{\alpha'}{2 \|\Omega\|_\omega} \Big) \nabla^2 (\Rm * \Rm)}{\nabla \Rm}^{\overline{j}} \cdot \Big[ \alpha' \Big( |\Rm|^2 + |D T|^2 \Big) + \mu \Big] \Bigg) \\
    &= 2 \alpha' p\, \Re \Bigg( \int_X \overline{\nabla}_{\overline{j}} \binner{G^{p-1} \cdot \Big( \frac{1}{2 \|\Omega\|_\omega} \Big) \cdot \Big[ \alpha' \Big( |\Rm|^2 + |D T|^2 \Big) + \mu \Big] \cdot \nabla^2 (\Rm * \Rm)}{\nabla \Rm}^{\overline{j}} \Bigg) \\*
    &\qquad - 2 \alpha' p\, \Bigg( \int_X \overline{\nabla} G^{p-1} \cdot \binner{\Big( \frac{1}{2 \|\Omega\|_\omega} \Big) \cdot \Big[ \alpha' \Big( |\Rm|^2 + |D T|^2 \Big) + \mu \Big] \cdot \nabla^2 (\Rm * \Rm)}{\nabla \Rm} \Bigg) \\*
    &\qquad - 2 (\alpha')^2 p\, \Bigg( \int_X G^{p-1} \cdot \binner{\Big( \frac{1}{2 \|\Omega\|_\omega} \Big) \cdot \nabla^2 (\Rm * \Rm)}{\nabla \Rm} \cdot \overline{\nabla}\Big( |\Rm|^2 + |D T|^2 \Big) \Bigg) \\
    &\leq C \int_X G^{p-1} \cdot \Big[ \alpha' \Big( |\Rm|^2 + |D T|^2 \Big) + \mu \Big] \cdot \Big( |\Rm| \cdot |D^2 \Rm| + |D \Rm|^2 \Big) \cdot |D \Rm|  \\*
    &\qquad + 2 a_0 B \alpha' p (p-1) \int_X G^{p-2} \cdot |\overline{\nabla} G| \cdot \Big[ \alpha' \Big( |\Rm|^2 + |D T|^2 \Big) + \mu \Big] \cdot \Big( 2 |\Rm| \cdot |D^2 \Rm| + 2 |D \Rm|^2 \Big) \cdot |D \Rm| \\*
    &\qquad + 2 a_0 B (\alpha')^2 p \int_X G^{p-1} \cdot \Big( 2 |\Rm| \cdot |D^2 \Rm| + 2 |D \Rm|^2 \Big) \cdot |D \Rm| \\*
    &\qquad \qquad \qquad \qquad \qquad \cdot \Big( C |\Rm| \cdot |DT| + 2 |\Rm| \cdot |D \Rm| + 4 |D T| \cdot |D^2 T| \Big)
\end{align*}
} A similar computation to the $\mathbf{(A)}$ term then gives that
{\allowdisplaybreaks
\begin{align*}
    &2p\, \Re \Bigg( \int_X G^{p-1} \cdot \overline{\nabla}_{\overline{j}} \binner{\Big( \frac{\alpha'}{2 \|\Omega\|_\omega} \Big) \nabla^2 (\Rm * \Rm)}{\nabla \Rm}^{\overline{j}} \cdot \Big[ \alpha' \Big( |\Rm|^2 + |D T|^2 \Big) + \mu \Big] \Bigg) \\
    &\leq C \eps^{-1} \int_X G^{p-1} + C \eps^{-1} \int_X G^p \\*
    &\qquad + C \eps \int_X G^{p-1} \cdot \Big( |D \Rm|^2 + |D^2 T|^2 \Big)^2 + C \eps \int_X G^{p-1} \cdot \Big( |D^2 \Rm|^2 + |D^3 T|^2 \Big) \\*
    &\qquad + 2 a_0 B \alpha' p (p-1) \int_X G^{p-3} \cdot \Big[ \alpha' \Big( |\Rm|^2 + |D T|^2 \Big) + \mu \Big] \cdot |D \Rm|^2 \cdot |\overline{\nabla} G|^2 \\*
    &\qquad + 2 a_0 B \alpha' p (p-1) \int_X G^{p-3} \cdot G^2 \cdot  \Big[ \alpha' \Big( |\Rm|^2 + |D T|^2 \Big) + \mu \Big] \cdot \Big( |D \Rm|^2 + |D^2 T|^2 \Big)^2 \\*
    &\qquad + 2 a_0 B C_0^2 \alpha' p (p-1) \int_X G^{p-3} \cdot G^2 \cdot  \Big[ \alpha' \Big( |\Rm|^2 + |D T|^2 \Big) + \mu \Big] \cdot \Big( |D^2 \Rm|^2 + |D^3 T|^2 \Big) \\*
    &\qquad + 12 a_0 B C_0^2 (\alpha')^2 p \int_X G^{p-1} \cdot \Big( |D \Rm|^2 + |D^2 T|^2 \Big)^2 + 12 a_0 B C_0^2 (\alpha')^2 p \int_X G^{p-1} \cdot \Big( |D^2 \Rm|^2 + |D^3 T|^2 \Big) \\*
    &\qquad + 16 a_0 B C_0 (\alpha')^2 p \int_X G^{p-1} \cdot \Big( |D \Rm|^2 + |D^2 T|^2 \Big)^2 \\
    &\leq C \eps^{-1} \int_X G^{p-1} + C \eps^{-1} \int_X G^p \\*
    &\qquad + 2 a_0 B \alpha' p (p-1) \int_X G^{p-2} \cdot |\overline{\nabla} G|^2 \\*
    &\qquad + \Big[ C \eps + 4 a_0 B C_0^2 (\alpha')^2 p (p-1) \\*
    &\qquad \qquad + 12 a_0 B C_0^2 (\alpha')^2 p + 16 a_0 B C_0 (\alpha')^2 p + 2 a_0 B \alpha' \mu p (p-1) \Big] \int_X G^{p-1} \cdot \Big( |D \Rm|^2 + |D^2 T|^2 \Big)^2 \\*
    &\qquad + \Big[ C \eps + 4 a_0 B C_0^4 (\alpha')^2 p (p-1) \\*
    &\qquad \qquad + 12 a_0 B C_0^2 (\alpha')^2 p + 2 a_0 B C_0^2 \alpha' \mu p (p-1) \Big] \int_X G^{p-1} \cdot \Big( |D^2 \Rm|^2 + |D^3 T|^2 \Big). \numberthis
\end{align*}
}Here we have again used that $\alpha' \Big( |\Rm|^2 + |DT|^2 \Big) + \mu \leq 2C_0^2 + \mu$.

The other parts of the $\mathbf{(B)}$ term have similar bounds and so we have
\begin{equation}
\begin{aligned}[b]
    &2p\, \Re \Bigg( \int_X G^{p-1} \cdot \mathbf{(B)} \Bigg) \\
    &\leq C \eps^{-1} \int_X G^{p-1} \cdot \Big( 1 + G \Big) + 10 a_0 B \alpha' p (p-1) \int_X G^{p-2} \cdot |\overline{\nabla} G|^2 \\
    &\qquad + \Big[ C \eps + 20 a_0 B C_0^2 (\alpha')^2 p (p-1) + 60 a_0 B C_0^2 (\alpha')^2 p \\
    &\qquad \qquad + 80 a_0 B C_0 (\alpha')^2 p + 10 a_0 B \alpha' \mu p (p-1) \Big] \int_X G^{p-1} \cdot \Big( |D \Rm|^2 + |D^2 T|^2 \Big)^2 \\
    &\qquad + \Big[ C \eps + 20 a_0 B C_0^4 (\alpha')^2 p (p-1) \\
    &\qquad \qquad + 60 a_0 B C_0^2 (\alpha')^2 p + 10 a_0 B C_0^2 \alpha' \mu p (p-1) \Big] \int_X G^{p-1} \cdot \Big( |D^2 \Rm|^2 + |D^3 T|^2 \Big).
\end{aligned}
\end{equation}

Combining what we have from the above, we get
\begin{equation}
\begin{aligned}[b]
    \del_t \Bigg( \int_X G^p \Bigg) &\leq C \eps^{-1} \int_X G^{p-1} \cdot \Big( 1 + G \Big) \\
    &\qquad + \Big[ C \eps + 6 a_0 B C_0^2 (\alpha')^2 \mu^{-1} p (p-1) + 10 a_0 B \alpha' p (p-1) - B^{-1} p (p-1) \Big] \int_X G^{p-2} \cdot |\overline{\nabla} G|^2 \\
    &\qquad + \Big[ C \eps + 140 a_0 B C_0^2 (\alpha')^2 p + 26 a_0 B C_0^2 (\alpha')^2 p (p-1) + 86 a_0 B C_0 (\alpha')^2 p \\
    &\qquad \qquad + 10 a_0 B \alpha' \mu p^2 + \frac{1}{2} B^{-1} \alpha' p - B^{-1} \alpha' p \Big] \int_X G^{p-1} \cdot \Big( |D \Rm|^2 + |D^2 T|^2 \Big)^2 \\
    &\qquad + \Big[ C \eps + 20 a_0 B C_0^4 (\alpha')^2 p (p-1) + 8320 B^3 C_0^2 \alpha' p + 12 a_0 B C_0^3 (\alpha')^2 p \\
    &\qquad \qquad + 128 a_0 B C_0^2 (\alpha')^2 p + 10 a_0 B C_0^2 \alpha' \mu p (p-1) \\
    &\qquad \qquad + 6 a_0 B C_0 \alpha' \mu p + 4 a_0 B \alpha' \mu p - B^{-1} \mu p \Big] \int_X G^{p-1} \cdot \Big( |D^2 \Rm|^2 + |D^3 T|^2 \Big).
\end{aligned}
\end{equation}

Thus, reading off the coefficients, we see that if
\begin{equation}
\label{eqn-alpha'-mu-1}
    6 a_0 B C_0^2 (\alpha')^2 \mu^{-1} + 10 a_0 B \alpha' < B^{-1},
\end{equation}
\begin{equation}
\label{eqn-alpha'-mu-2}
\begin{aligned}[b]
    &140 a_0 B C_0^2 \alpha' + 26 a_0 B C_0^2 \alpha' (p-1) + 86 a_0 B C_0 \alpha' + 10 a_0 B \mu p + \frac{1}{2} B^{-1} < B^{-1},
\end{aligned}
\end{equation}
\begin{equation}
\label{eqn-alpha'-mu-3}
\begin{aligned}[b]
    &20 a_0 B C_0^4 (\alpha')^2 (p-1) + 8320 B^3 C_0^2 \alpha' + 12 a_0 B C_0^3 (\alpha')^2 \\
    &\qquad + 128 a_0 B C_0^2 (\alpha')^2 + 10 a_0 B C_0^2 \alpha' \mu (p-1) + 8 a_0 B C_0 \alpha' \mu + 4 a_0 B \alpha' \mu < B^{-1} \mu,
\end{aligned}
\end{equation}
then we can absorb the corresponding terms into the negative ones.

We see that if
\begin{equation}
    \alpha' = O \Big( a_0^{-1} B^{-6} \max(1,C_0)^{-2} p^{-1} \Big) \text{ and } \mu = O \Big( a_0^{-1} B^{-2} p^{-1} \Big),
\end{equation}
then the three inequalities could hold.

To this end, we get that if
\begin{equation}
    \alpha' < \frac{1}{10^6 a_0 B^6 \max(1,C_0)^2 p},
\end{equation}
then by picking
\begin{equation}
    \mu = \frac{1}{100 a_0 B^2 p},
\end{equation}
we satisfy the inequalities \eqref{eqn-alpha'-mu-1} - \eqref{eqn-alpha'-mu-3} and can then choose $\eps = \eps(k,\alpha',p,\mu) > 0$ such that
\begin{equation}
    \del_t \Bigg( \int_X G^p \Bigg) \leq C \int_X G^{p-1} + C \int_X G^p \leq C \int_X 1 + C \int_X G^p.
\end{equation}

As a consequence of Gr\"{o}nwall's inequality again, we get
\begin{thm}
\label{thm-G^p}
    Let $k = 1$ and set $G_r = |D^r \Rm|^2 + |D^{r+1} T|^2$. Given the assumptions \eqref{eqn-assumptions-1} - \eqref{eqn-assumptions-3}, define
    \begin{equation}
        G = \Big[ \alpha' \Big( |\Rm|^2 + |DT|^2 \Big) + \mu \Big] \cdot \Big( |D \Rm|^2 + |D^2 T|^2 \Big) = (\alpha' G_0 + \mu) \cdot G_1,
    \end{equation}
    for some $\mu$. If:

    \begin{enumerate}[label=(\alph*)]
        \item
        $p \geq 3$ and
        \begin{equation}
            \alpha' < \frac{1}{10^6 a_0 B^6 \max(1,C_0)^2 p},
        \end{equation}
        then for $\mu = \frac{1}{100 a_0 B^2 p}$ there exists a constant $\Lambda_p = \Lambda_p(k, \alpha') > 0$ such that
        \begin{equation}
        \label{eqn-G^p}
            \int_X G^p(t) \leq \Bigg( 1 + \int_X G^p(0) \Bigg) e^{\Lambda_p t} < \Bigg( 1 + \int_X G^p(0) \Bigg) e^{\Lambda_p \tau}.
        \end{equation}
        That is, $\int_X G^p(t)$ is uniformly bounded along the flow.
    
        In particular, we get that
        \begin{equation}
            \int_M G_1^p(t) \text{ and both } \Bigg( \int_X |D \Rm(t)|^{2p} \Bigg)^{\frac{1}{2p}} \text{ and } \Bigg( \int_X |D^2 T(t)|^{2p} \Bigg)^{\frac{1}{2p}}
        \end{equation}
        are bounded along the anomaly flow;

        \item
        instead $1 \leq p < 3$ and 
        \begin{equation}
            \alpha' < \frac{1}{3 \cdot 10^6 a_0 B^6 C_0^2} \text{ }\Bigg(\text{or } \alpha' < \frac{1}{3 \cdot 10^6 a_0 B^6} \text{ if } C_0 < 1 \Bigg),
        \end{equation}
        then for $\mu = \frac{1}{300 a_0 B^2}$, the function $\int_X G^p(t)$ and also
        \begin{equation}
            \int_X G_1^p(t) \text{ and both }
            \Bigg( \int_X |D \Rm(t)|^{2p} \Bigg)^{\frac{1}{2p}} \text{ and } \Bigg( \int_X |D^2 T(t)|^{2p} \Bigg)^{\frac{1}{2p}}
        \end{equation}
        are uniformly bounded along the anomaly flow.
    \end{enumerate}
    \blktr
\end{thm}

\subsubsection{Estimates on \texorpdfstring{$|D^{k+1} \Rm|$}{|D(k+1) Rm|} and \texorpdfstring{$|D^{k+2} T|$}{|D(k+2) T|}}
\label{subsubsect-higher-order-ests-k=1}

One can check that
\begin{equation}
\begin{aligned}[b]
    \del_t \Big( |D^2 \Rm|^2 + |D^3 T|^2 \Big) &\leq \frac{1}{2} \Big( \frac{1}{2 \|\Omega\|_\omega} \Big) \Delta_R \Big( |D^2 \Rm|^2 + |D^3 T|^2 \Big) - B^{-1} \Big( |D^3 \Rm|^2 + |D^4 T|^2 \Big) \\
    &\qquad + \sum_{m+l = 2} 2 \Re \Big( \overline{\nabla}_{\overline{j}} \binner{\Big( \frac{\alpha'}{2 \|\Omega\|_\omega} \Big) \nabla^{m+1} \overline{\nabla}^l (\Rm * \Rm)}{\nabla^m \overline{\nabla}^l \Rm}^{\overline{j}} \Big) \\
    &\qquad + \sum_{m' + l'= 3} 2 \Re \Big( \nabla_i \binner{\Big( \frac{\alpha'}{2 \|\Omega\|_\omega} \Big)  \nabla^{m'} \overline{\nabla}^{l'} (\Rm * \Rm)}{\nabla^{m'} \overline{\nabla}^{l'} T}^i \Big) \\
    &\qquad + C (\eps')^{-1} 1 + C (\eps')^{-1} \Big( |D \Rm|^2 + |D^2 T|^2 \Big) + C (\eps')^{-1} \Big( |D \Rm|^2 + |D^2 T|^2 \Big)^2 \\
    &\qquad + C (\eps')^{-1} \Big( |D^2 \Rm|^2 + |D^3 T|^2 \Big) \\
    &\qquad + \Big[ C \eps' + 42 a_0 B \alpha' \Big] \Big( |D \Rm|^2 + |D^2 T|^2 \Big) \cdot \Big( |D^2 \Rm|^2 + |D^3 T|^2 \Big) \\
    &\qquad + \Big[ C \eps' + 6 a_0 B C_0 \alpha' + 4 a_0 B \alpha' \Big] \Big( |D^3 \Rm|^2 + |D^4 T|^2 \Big).
\end{aligned}
\end{equation}

We choose a similar but more complicated test function this time and set
\begin{equation}
\begin{aligned}[b]
    G' = \Big[ \alpha' \Big( |\Rm|^2 + |D T|^2 \Big) + \mu' \Big] \cdot \Big( |D^2 \Rm|^2 + |D^3 T|^2 \Big) = (\alpha' G_0 + \mu') \cdot G_2.
\end{aligned}
\end{equation}
where $\mu'$ is a constant to be determined later. Analogous computations show that for $0 < \eps' < 1$,
\begin{equation}
\begin{aligned}[b]
    \del_t G' &\leq \frac{1}{2} \Big( \frac{1}{2 \|\Omega\|_\omega} \Big) (\Delta_R G') - B^{-1} \alpha' \Big( |D \Rm|^2 + |D^2 T|^2 \Big) \cdot \Big( |D^2 \Rm|^2 + |D^3 T|^2 \Big) \\
    &\qquad - B^{-1} \mu' \Big( |D^3 \Rm|^2 + |D^4 T|^2 \Big) \\
    &\qquad + \Big( \frac{\alpha'}{2 \|\Omega\|_\omega} \Big) 2 \Re \binner{\nabla \Big( |\Rm|^2 + |D T|^2 \Big)}{\nabla \Big( |D^2 \Rm|^2 + |D^3 T|^2 \Big)} \\
    &\qquad + 2 \Re \mathbf{(A')} + 2 \Re \mathbf{(B')} \\
    &\qquad + C (\eps')^{-1} \Big( 1 + G_1 + G_1^2 + G' \Big) \\
    &\qquad + \Big[ C \eps' + 6 a_0 B C_0 (\alpha')^2 + 84 a_0 B C_0^2 (\alpha')^2 + 42 a_0 B \alpha' \mu' \Big] \Big( |D \Rm|^2 + |D^2 T|^2 \Big) \cdot \Big( |D^2 \Rm|^2 + |D^3 T|^2 \Big) \\
    &\qquad + \Big[ C \eps' + 12 a_0 B C_0^3 (\alpha')^2 + 8 a_0 B C_0^2 (\alpha')^2 + 6 a_0 B C_0 \alpha' \mu' + 4 a_0 B \alpha' \mu' \Big] \Big( |D^3 \Rm|^2 + |D^4 T|^2 \Big),
\end{aligned}
\end{equation}
where the terms $\mathbf{(A')}$ and $\mathbf{(B')}$ are given by

\begin{equation}
\begin{aligned}[b]
    \mathbf{(A')} &= \overline{\nabla}_{\overline{j}} \binner{\Big( \frac{\alpha'}{2 \|\Omega\|_\omega} \Big) \nabla (\Rm * \Rm)}{\Rm}^{\overline{j}} \cdot \alpha' \Big( |D^2 \Rm|^2 + |D^3 T|^2 \Big) \\
    &\qquad + \sum_{m' + l'= 1} \nabla_i \binner{\Big( \frac{\alpha'}{2 \|\Omega\|_\omega} \Big)  \nabla^{m'} \overline{\nabla}^{l'} (\Rm * \Rm)}{\nabla^{m'} \overline{\nabla}^{l'} T}^i \cdot \alpha' \Big( |D^2 \Rm|^2 + |D^3 T|^2 \Big),
\end{aligned}
\end{equation}
\begin{equation}
\begin{aligned}[b]
    \mathbf{(B')} &= \sum_{m+l = 2} \overline{\nabla}_{\overline{j}} \binner{\Big( \frac{\alpha'}{2 \|\Omega\|_\omega} \Big) \nabla^{m+1} \overline{\nabla}^l (\Rm * \Rm)}{\nabla^m \overline{\nabla}^l \Rm}^{\overline{j}} \cdot \Big[ \alpha' \Big( |\Rm|^2 + |D T|^2 \Big) + \mu' \Big] \\
    &\qquad + \sum_{m' + l'= 3} \nabla_i \binner{\Big( \frac{\alpha'}{2 \|\Omega\|_\omega} \Big)  \nabla^{m'} \overline{\nabla}^{l'} (\Rm * \Rm)}{\nabla^{m'} \overline{\nabla}^{l'} T}^i \cdot \Big[ \alpha' \Big( |\Rm|^2 + |D T|^2 \Big) + \mu' \Big].
\end{aligned}
\end{equation}

Proceeding as before, we get that for $p \geq 3$,
\begin{equation}
\begin{aligned}[b]
    \del_t \Bigg( \int_X (G')^p \Bigg) &\leq C (\eps')^{-1} \int_X (G')^{p-1} + C (\eps')^{-1} \int_X (G')^{p-1} \cdot G_1 + C (\eps')^{-1} \int_X (G')^{p-1} \cdot G_1^2 + C \eps^{-1} \int_X (G')^p \\
    &\qquad + \Big[ C \eps' + 6 a_0 B C_0^2 (\alpha')^2 (\mu')^{-1} p (p-1) + 56 a_0 B \alpha' p (p-1) - B^{-1} p (p-1) \Big] \int_X G^{p-2} \cdot |\overline{\nabla} G|^2 \\
    &\qquad + \Big[ C \eps' + 252 a_0 B C_0^2 (\alpha')^2 p + 90 a_0 B C_0^2 (\alpha')^2 p (p-1) \\
    &\qquad \qquad + 510 a_0 B C_0 (\alpha')^2 p + 42 a_0 B \alpha' \mu' p^2 \\
    &\qquad \qquad + \frac{1}{2} B^{-1} \alpha' p - B^{-1} \alpha' p \Big] \int_X G^{p-1} \cdot \Big( |D \Rm|^2 + |D^2 T|^2 \Big) \cdot \Big( |D^2 \Rm|^2 + |D^3 T|^2 \Big) \\
    &\qquad + \Big[ C \eps' + 28 a_0 B C_0^4 (\alpha')^2 p (p-1) + 16000 B^3 C_0^2 \alpha' p + 12 a_0 B C_0^3 (\alpha')^2 p \\
    &\qquad \qquad + 176 a_0 B C_0^2 (\alpha')^2 p + 14 a_0 B C_0^2 \alpha' \mu' p (p-1) \\
    &\qquad \qquad + 6 a_0 B C_0 \alpha' \mu' p + 4 a_0 B \alpha' \mu' p - B^{-1} \mu' p \Big] \int_X G^{p-1} \cdot \Big( |D^2 \Rm|^2 + |D^3 T|^2 \Big).
\end{aligned}
\end{equation}

As in \S \ref{subsect-k=2}, we also require $\int_X G_1^{2p}$ to be uniformly bounded and must have $\alpha'$ satisfy an improved bound from Theorem \ref{thm-G^p}. Similar reasoning then yields
\begin{cor}
\label{cor-G'^p}
    Let $k = 1$ and set $G_r = |D^r \Rm|^2 + |D^{r+1} T|^2$. Given the assumptions \eqref{eqn-assumptions-1} - \eqref{eqn-assumptions-3}, define
    \begin{equation}
        G' = \Big[ \alpha' \Big( |\Rm|^2 + |DT|^2 \Big) + \mu' \Big] \cdot \Big( |D^2 \Rm|^2 + |D^3 T|^2 \Big) = (\alpha' G_0 + \mu') \cdot G_2,
    \end{equation}
    for some $\mu'$. If:

    \begin{enumerate}[label=(\alph*)]
        \item
        $p \geq 3$ and
        \begin{equation}
            \alpha' < \frac{1}{10^7 a_0 B^6 \max(1,C_0)^2 p},
        \end{equation}
        then for $\mu' = \frac{1}{100 a_0 B^2 p}$ there exists a constant $\Lambda_p' = \Lambda_p'(k, \alpha') > 0$ such that
        \begin{equation}
        \label{eqn-G'^p}
            \int_X (G')^p(t) \leq \Bigg( 1 + \int_X (G')^p(0) \Bigg) e^{\Lambda_p' t} < \Bigg( 1 + \int_X (G')^p(0) \Bigg) e^{\Lambda_p' \tau}.
        \end{equation}
        That is, $\int_X G^p(t)$ is uniformly bounded along the flow.
    
        In particular, we get that
        \begin{equation}
            \int_M G_2^p(t) \text{ and both } \Bigg( \int_X |D^2 \Rm(t)|^{2p} \Bigg)^{\frac{1}{2p}} \text{ and } \Bigg( \int_X |D^3 T(t)|^{2p} \Bigg)^{\frac{1}{2p}}
        \end{equation}
        are bounded along the anomaly flow;

        \item
        instead $1 \leq p < 3$ and 
        \begin{equation}
            \alpha' < \frac{1}{3 \cdot 10^7 a_0 B^6 \max(1,C_0)^2},
        \end{equation}
        then for $\mu' = \frac{1}{300 a_0 B^2}$, the function $\int_X (G')^p(t)$ and also
        \begin{equation}
            \int_X G_2^p(t) \text{ and both }
            \Bigg( \int_X |D^2 \Rm(t)|^{2p} \Bigg)^{\frac{1}{2p}} \text{ and } \Bigg( \int_X |D^3 T(t)|^{2p} \Bigg)^{\frac{1}{2p}}
        \end{equation}
        are uniformly bounded along the anomaly flow.
    \end{enumerate}
    \blktr
\end{cor}

Using the Sobolev embedding theorem and inductive bootstrapping argument from \S \ref{subsect-Sobolev-induction}, we can conclude the following:

\begin{thm}
\label{thm-L^infty-k=1}
    Suppose that there exist positive constants $B, C_0$ such that
    \begin{equation}
        B^{-1} \leq \Big( \frac{1}{2 \|\Omega\|_\omega} \Big) \leq B,
    \end{equation}
    \begin{equation}
        |T|, |\overline{T}|, |\Rm|, |DT|, |D \overline{T}| \leq C_0
    \end{equation}
    along the anomaly flow on $t \in [0, \tau)$. If 
    \begin{equation}
            \alpha' < \frac{1}{3 \cdot 10^7 a_0 B^6 \max(1,C_0)^2},
        \end{equation}
    then there exist positive constants $C_q$ for $q \geq 1$ such that
    \begin{equation}
        |D^q \Rm|, |D^{q+1} T|, |D^{q+1} \overline{T}| \leq C_q
    \end{equation}
    along the anomaly flow on $t \in [0,\tau)$. These constants $C_q$ depend on the initial metric $g_0$, the slope parameter $\alpha'$, and the initial bounds $B, C_0$.
    \blktr
\end{thm}

\begin{rmk}
\label{rmk-dimensionality}
    Instead of our choice of coupled upper and lower bounds
    \begin{equation}
        B^{-1} \leq \Big( \frac{1}{2 \|\Omega\|_\omega} \Big) \leq B,
    \end{equation}
    we could have instead chosen independent bounds
    \begin{equation}
        B_{min} \leq \Big( \frac{1}{2 \|\Omega\|_\omega} \Big) \leq B_{max}.
    \end{equation}

    In this case, our derived bounds on $\alpha'$ from Theorem \ref{thm-L^infty-k=1} (the $k = 1$ case) and Corollary \ref{cor-L^infty-k=2} (the $k \geq 2$ case) respectively become
    \begin{equation}
         \alpha' < \frac{\min(1,B_{min})^3}{3 \cdot 10^7 a_0 \max(1,B_{max})^3 \max(1,C_0)^2}, \qquad (k = 1),
    \end{equation}
    \begin{equation}
        \alpha' < \frac{B_{min}}{26 a_0 B_{max} C_0}, \qquad (k \geq 2).
    \end{equation}
    
    In particular these derived conditions can respectively be rewritten as
    \begin{equation}
        \alpha' \cdot |\Rm|^2 < \Pi_1, \qquad (k = 1),
    \end{equation}
    \begin{equation}
        \alpha' \cdot |\Rm| < \Pi_2, \qquad (k \geq 2),
    \end{equation}
    for some dimensionless constants $\Pi_1$ and $\Pi_2$. This is intriguing as the units on the LHS of the previous two equations differ but the RHS of both are dimensionless.
    \blktr
\end{rmk}

\section{Long-time Existence}
\label{sect-long-time-existence}

Now that we have established the relevant $L^\infty$-bounds for the covariant derivatives and torsion, we can appeal to the argument in \cite{PPZ18b} to obtain long-time existence of the anomaly flow. For completeness, we outline the argument below.

\subsection{\texorpdfstring{$C^0$}{C0}- and \texorpdfstring{$C^1$}{C1}-Bounds}
\label{subsect-C0-C1-bds}

Suppose that there exist positive constants $B$, $C_0$ such that
\begin{equation}
    B^{-1} \leq \Big( \frac{1}{2 \|\Omega\|_\omega} \Big) \leq B,
\end{equation}
\begin{equation}
    |T|, |\overline{T}|, |\Rm|, |DT|, |D \overline{T}| \leq C_0
\end{equation}
along the anomaly flow on $[0,\tau)$. The metrics $g$ are uniformly equivalent close to time $\tau$, and so if $\what{g}$ is a fixed reference metric, the endomorphism
\begin{equation}
    {h^\alpha}_\beta = (\what{g})^{\alpha \overline{\gamma}} g_{\overline{\gamma} \beta}
\end{equation}
has a uniform $C^0$-bound.

We have that the Chern curvatures of $\what{g}$ and $g$ are related by
\begin{equation}
    \tensor{\what{R}}{_{\overline{p}}_q^\alpha_\beta} - \tensor{R}{_{\overline{p}}_q^\alpha_\beta} = \del_{\overline{p}} ({h^\alpha}_\gamma \what{\nabla}_q {h^\gamma}_\beta).
\end{equation}
Since $\Rm$ is bounded, we have a $C^1$-bound on the endomorphism $h$.

\subsection{\texorpdfstring{$C^k$}{Ck}-Bounds}
\label{subsect-Ck-bds}

For each $q$ and for a tensor $A$, set
\begin{equation}
    |\what{D}^q A|^2 = \sum_{m+l = q} |\what{\nabla}^m \overline{\what{\nabla}^l} A|^2.
\end{equation}
We also set the tensor $S$ to be the difference of the Christoffel symbols of the reference and moving metrics.
\begin{equation}
    \tensor{S}{^\alpha_k_\beta} = \tensor{\what{\Gamma}}{^\alpha_k_\beta} - \tensor{\Gamma}{^\alpha_k_\beta} = - g^{\alpha \overline{\gamma}} \what{\nabla}_k g_{\overline{\gamma} \beta}.
\end{equation}

From Theorem \ref{thm-L^infty-k=1}, we see that if $\alpha'$ is sufficiently small, then all covariant derivatives of $\Rm$ and $T$ with respect to the evolving metric are uniformly bounded on $[0,\tau)$. Since
\begin{equation}
    \del_t g = \frac{1}{2 \|\Omega\|_\omega} \Bigg[ \Rm + T * \overline{T} + \alpha' \Big( \Rm * \Rm + \Phi \Big) \Bigg],
\end{equation}
it follows that all covariant derivatives of $\del_t g$ with respect to the evolving metric are also uniformly bounded.

We also have the following from \cite{PPZ18b},
\begin{propn}[\cite{PPZ18b} Proposition 2.]
\label{propn-D^kS-hatD^k+1g}
    Suppose all covariant derivatives of $\Rm$ and $T$ with respect to the evolving metric are uniformly bounded along the anomaly flow on $[0,\tau)$. If for $k \geq 1$ there exist positive constants $C_0',C_1',\ldots,C_q'$ such that
    \begin{equation}
        |D^q S|, |\what{D}^{q+1} g| \leq C_q' \text{ for } 1 \leq q \leq k-1,
    \end{equation}
    \begin{equation}
        |g|, |S|, |\hat{D} g| \leq C_0'
    \end{equation}
    along the anomaly flow, then there exists some positive $C_k'$ such that
    \begin{equation}
        |D^k S|, |\what{D}^{k+1} g| \leq C_k'.
    \end{equation}
    \blktr
\end{propn}

The $C^0$- and $C^1$- bounds from \S \ref{subsect-C0-C1-bds} imply the existence of some positive $C_0'$ such that
\begin{equation}
    |g|, |S|, |\what{D} g| \leq C_0'. 
\end{equation} 
Inductively applying Proposition \ref{propn-D^kS-hatD^k+1g}, we see that the covariant derivatives of $g$ with respect to the reference metric $\what{g}$, that is $|\what{D}^q g|$ for each $q \geq 0$, are all uniformly bounded along the flow.  

Now, for any $i$, $m$, $l$, we have
\begin{equation}
    \del_t^i (\what{\nabla}^m \overline{\what{\nabla}^l} g) = \what{\nabla}^m \overline{\what{\nabla}^l} \del_t^i g = \what{\nabla}^m \overline{\what{\nabla}^l} \del_t^{i-1} \Bigg[ \Big( \frac{1}{2 \|\Omega\|_\omega} \Big) \Big[ \Rm + T * \overline{T} + \alpha' (\Rm * \Rm + \Phi) \Big] \Bigg]. 
\end{equation}
The time derivatives of $\del_t g$ can also be expressed in terms of time derivatives of connections, $\Rm$ and $T$, which our calculations in previous sections have shown to be bounded. As such
\begin{equation}
    \del_t^i (\what{\nabla}^m \overline{\what{\nabla}^l} g)
\end{equation}
is uniformly bounded along the flow on $[0,\tau)$. We can thus extend the anomaly flow smoothly to the time $t = \tau$ and by the short-time existence of the flow from Theorem 2 of \cite{PPZ18c}, we can further extend it to $[0,\tau + \eps)$ for some $\eps > 0$.

We have thus shown the following:
\begin{thm}
\label{thm-main}
    Suppose that there exist positive constants $B, C_0$ such that
    \begin{equation}
        B^{-1} \leq \Big( \frac{1}{2 \|\Omega\|_\omega} \Big) \leq B,
    \end{equation}
    \begin{equation}
        |T|, |\overline{T}|, |\Rm|, |DT|, |D \overline{T}| \leq C_0
    \end{equation}
    along the anomaly flow on $t \in [0, \tau)$. If 
    \begin{equation}
            \alpha' < \frac{1}{3 \cdot 10^7 a_0 B^6 \max(1,C_0)^2},
        \end{equation}
    then the flow can be extended to $[0, \tau + \eps)$ for some $\eps > 0$.
    \blktr
\end{thm}

\appendix

\section{Identities for Hermitian Metrics and Chern Connections}
\label{sect-useful-ids}

In this appendix, we list some useful identities that will be used often.

The conformally balanced condition $d(\|\Omega\|_\omega \omega^2) = 0$ is equivalent to
\begin{equation}
    T_i = \nabla_i \log \|\Omega\|_\omega, \qquad \overline{T}_{\overline{j}} = \overline{\nabla}_{\overline{j}} \log \|\Omega\|_\omega.
\end{equation}
As such, we see that
\begin{equation}
\label{eqn-nabla-1/2-Omega}
\begin{aligned}[b]
    \nabla \Big( \frac{1}{2 \|\Omega\|_\omega} \Big) = - \Big( \frac{1}{2 \|\Omega\|_\omega} \Big) \cdot T, \qquad \overline{\nabla} \Big( \frac{1}{2 \|\Omega\|_\omega} \Big) = - \Big( \frac{1}{2 \|\Omega\|_\omega} \Big) \cdot \overline{T}.
\end{aligned}
\end{equation}

Repeated application of the above yields the following:
\begin{equation}
\label{eqn-nabla^m-nabla-bar^l-1/2-Omega}
\begin{aligned}[b]
    \nabla^m \overline{\nabla}^l \Big( \frac{1}{2 \|\Omega\|_\omega} \Big) = \Big( \frac{1}{2 \|\Omega\|_\omega} \Big) \sum_{\substack{i_1 + \ldots + i_r + (r-s) = m \\ j_1 + \ldots + j_s + s = l}} \nabla^{i_1} \overline{\nabla}^{j_1} \overline{T} * \ldots *  \nabla^{i_s} \overline{\nabla}^{j_s} \overline{T} * \nabla^{i_{s+1}} T * \ldots * \nabla^{i_r} T.
\end{aligned}
\end{equation}

We also note the general commutator identities: For a generic tensor $A$,
\begin{equation}
\label{eqn-nabla^m-nabla-bar^l-Laplace}
\begin{aligned}[b]
    \nabla^m \overline{\nabla}^l (\Delta_R A) &= \Delta_R (\nabla^m \overline{\nabla}^l A) + \sum_{i=0}^m \sum_{j=0}^l (\nabla^{m-i} \overline{\nabla}^{l-j} A) * (\nabla^i \overline{\nabla}^j \Rm) \\
    &\qquad + \sum_{i=0}^m \sum_{j=0}^l (\nabla^{m-i} \overline{\nabla}^{l+1-j} A) * (\nabla^i \overline{\nabla}^j T) + \sum_{i=0}^m \sum_{j=0}^l (\nabla^{m+1-i} \overline{\nabla}^{l-j} A) * (\nabla^i \overline{\nabla}^j \overline{T}),
\end{aligned}
\end{equation}
and
\begin{equation}
\label{eqn-nabla-bar^l-nabla^m}
\begin{aligned}[b]
    \overline{\nabla}^l \nabla^m A &= \sum_{r=0}^{\min(m,l)} \sum_{\substack{i_0 + \ldots + i_r = m-r \\ j_0 + \ldots + j_r = l-r}} (\nabla^{i_0} \overline{\nabla}^{j_0} A) * (\nabla^{i_1} \overline{\nabla}^{j_1} \Rm) * (\nabla^{i_r} \overline{\nabla}^{j_r} \Rm).
\end{aligned}
\end{equation}

Lastly, we have the divergence theorem for the Chern connection
\begin{equation}
\label{eqn-div-thm}
    \int_X \nabla_i V^i = \int_X T_i \cdot V^i \text{ and } \int_X \overline{\nabla}_{\overline{j}} V^{\overline{j}} = \int_X \overline{T}_{\overline{j}} \cdot V^{\overline{j}}.
\end{equation}

\end{document}